\documentclass[12pt]{article}

\title{Diophantine Approximation on varieties III: Approximation of non-algebraic points by algebraic points}

\author{Heinrich Massold}

\usepackage{bezier,amsfonts}

\usepackage{amsmath,amsfonts,bbm,amssymb}

\newtheorem{Satz}{}[section]

\newcommand{\satz}[1]{\vspace{2mm} \begin{Satz}{\bf #1}}

\newcommand{\proof}{\vspace{3mm} {\sc Proof}\hspace{0.2cm}}

\newcommand{\la}{\langle}

\newcommand{\ra}{\rangle}

\newcommand{\di}{{\mbox{div}}}

\newcommand{\rk}{\mbox{rk}}

\newcommand{\R}{{\mathbbm R}}

\newcommand{\Pe}{{\mathbbm P}}

\newcommand{\Z}{\mathbbm{Z}}

\newcommand{\N}{\mathbbm{N}}

\newcommand{\CC}{{\cal C}}

\newcommand{\CH}{{\cal H}}

\newcommand{\CW}{{\cal W}}

\newcommand{\CO}{{\cal O}}

\newcommand{\CX}{{\cal X}}

\newcommand{\CY}{{\cal Y}}

\newcommand{\CZ}{{\cal Z}}

\newcommand{\CL}{{\cal L}}

\newcommand{\CF}{{\cal F}}

\newcommand{\Q}{{\mathbbm Q}}

\newcommand{\C}{{\mathbbm C}}

\newcommand{\wde}{\widehat{\deg}}

\newcommand{\codim}{\mbox{codim}}

\newcommand{\spec}{\mbox{Spec}}

\begin{document}

\parindent0mm

\maketitle

\thispagestyle{empty}

\tableofcontents

\section{Introduction}

Let $k/\Q$ be a number field with ring of integers $\CO_k$, and $\CX$ a flat
quasi projective scheme of finite type over $\spec \; \CO_k$, equipped
with an ample metrized line bundle $\bar{\CL}$.
The basse extension of $\CX$ to $k$ will be denoted by $X$. 
Let further $\sigma: k \to \C$ be some embeding, $X(\C_\sigma)$ the 
$\C_\sigma$- valued points of $\CX$ which usually will also be denoted by
$X$, and $|\cdot, \cdot|$ any metric on $\CX(\C_\sigma)$ that induces the
usual $\C$-topology.

\satz{Proposition} \label{alg}
In the above situation, for every $\beta \in \CX(\bar{k})$ with 
$\beta \neq \alpha$,
\[ \log |\alpha,\beta| \geq -c_1 \deg \beta  c_2 h(\beta)), \]
where $c_1$ is positive number depnding only on the height of $\alpha$, and
$c_2$ is a positive number depending only on the degree of $\beta$.
\end{Satz}

\vspace{2mm}
A proof that uses
Liouvilles Theorem, and the relation between distance and algebraic
distance in Theorem \ref{bezout}.1, from \cite{App1} is given in the appendix.

The Proposition says that the best approximation of a point 
$\alpha \in \CX(\bar{k})$ by 
another algebraic point $\beta$ is linear in the degree and height of $\beta$.
If on the other hand $\theta \in X(\C_\sigma)$, is not algebraic, a
possible approximation of $\theta$ by algebraic points is much better.
How good the approximation is,
depends mainly on the transecendence degree of the field generated by $\theta$.
More specifically there is the following Theorem whose proof is the main
objecture of this paper

\satz{Theorem} \label{main}
In the above situation, there exists a positive real number $b$, such
that for every sufficently big real number $a$, and
any nonalgebraic point $\theta \in \CX(\C_\sigma)$ with Zariski 
closure denoted by $\CY$, there exists  an infinite subset 
$M \subset \N$ such that for all $D \in M$ there is an
algebraic point $\alpha_D \in \CY(\bar{k})$ fullfilling
\[ \deg(\alpha_D) \leq D^t, \quad h(\alpha_D) \leq a D^t, \quad
   \log |\alpha_D, \theta| \leq -  
   \frac{1}{\left(\frac{h(\CX)}a + 2 \deg X\right)^{\frac{t+1}t}} \; b a 
   D^{t+1}, \]
where $t$ denotes the relative dimension of $\CY$ over $\spec \; \CO_k$.
\end{Satz}

\satz{Corollary} \label{coreinl}
For $\CX = \Pe^t$,
$D \in M$, and every effective cycle $\CY_D$ whoose support
contains the support of $\alpha_D$ we have
\[ D(X_D,\theta) \leq \log |X_D,\theta| \leq - b a D^{t+1} \]
\end{Satz}

\proof
1. Follows immediately from the trivial fact
$|X_D,\theta| \leq |\alpha_D,\theta|$ if $supp(\alpha_D) \subset supp(X_D)$, 
Theorem \ref{bezout}.1, and chaning $b$ slightly.



\vspace{2mm}

Of course it would be desireable to have an analogue for the lower
bound of the approximations that Proposition \ref{alg} supplies when
$t$ equals zero, i.\@ e.\@ when $\theta$ is algebraic.
thereby saying that the approximation given by 
Theorem \ref{main} up to a constant is best possible.
However such a lower estimate is not possible in 
general due to the existence of transcendental points like Liouville numbers.
The best that can be obtained is:

\satz{Conjecture} \label{con}
There is a constant $c$ depending only on $t$, such that for almost
all $\theta \in \CX(\C_\sigma)$ for sufficiently big $D$
the inequality $\log |\alpha_D,\theta| \leq - a b D^{t+1}$ implies
\[ \deg \alpha_D \geq c D^t \quad \mbox{or} \quad h(\alpha_D) \geq c a D^t. \]
\end{Satz}

A proof of this conjecture will be the subject of \cite{Mahler}.

\vspace{2mm}

There are several applications of this Theorem to algebraic independence
theory and other topics in transcendence theory, some of which are alluded
to in the outlook at the end of the paper.

\vspace{2mm}

The proof of Theorem \ref{main} has two steps. In the 
first step using the Theorem of Minkowski, estimates for the algebraic
and arithmetic Hilbert functions,
and the arithmetic Bezout Theorem from \cite{App1}, for each $D \in \N$,
cycles of codimension 
$t$ are construed that have small algebraic distance to $\theta$.
The main idea is to successively intersect properly intersecting cycles of 
small algebraic distance. If $\CY$ is such an effective cycle with
small algebraic distance, bounden height and degree, has codimension one,
and is irreducible,  a hypersurface properly intersecting $\CY_s$ can
be obtained if $\CY$ fullfills certain regularity conditions. Then, 
the intersection of $\CY$, and this hypersurface, by the metric B\'ezout
Theorem has also small algebraic distance to $\theta$.

If $\CY$ does not fullfill the regularity conditions, one nonetheless can find
a cycle $\CX$ that contains $\CY$ and has good regularity, and work with
this cycle. Again either there then is a cycle of codimension one in $\CY_s$ 
that or $\CX$ good approximation properties but small degree and height
with respect to a certain measure. Since the degree and height can not
infinitely decrease, one finds the disired cycle of codimension $t$.

The second step uses the fact that the algebraic distance of a cylcle $X$
not containing $\theta$ essentially equals the sum of logarithms of
certiain points on $X$ to $\theta$, and if $X$ is $0$-dimensional, this
is just the sum of logarithms of distances of the points in $X$ counted
with multiplicity to $\theta$. It is shown 
that for $D$ in an infinite set and $\CY_D$ a cycle of codimension $t$ with
good approximation properties with respect to $\theta$,
the sizes of these distances are very unequally distributed
to the effect that a fixed portion of the algebraic distance actually
comes from the logarithm of the distance of the $\C$-valued point of $\CY_D$
closest to $\theta$. A general version of the Liouville inequality plays an 
essential role in this step.

The starting point for the proof presented here was the proof of Theorem
\ref{main} in
\cite{RW} for $t=1$. Like that one it depends on a certain metric 
B\'ezout Theorem that relates distances and algebraic distances of
properly intersecting cycles in projective space with certain multiplicities
that are defined by the metric position of the cycles with respect to $\theta$.
This metric B\'ezout Theorem was proved in \cite{RW} for $t=1$, and 
generalized to arbitrary dimension in \cite{App1}. Some of the new problems
that appear in the higher dimensional case can also be solved using
these multiplicities.

The other significance difference from the one dimensional case, it the
fact that in higher dimensions for each $D>>0$, one needs to construe 
several in particular higher codimensional 
cycles with good approximation that fullfill certain conditions
of proper intersection.
The key ingredient to attain this, is the
concept of locally complete intersections in projective space
as introduced and investigated
in \cite{CP} and applied to Diophantine Approximation in \cite{Ph}.
Only with them it is possible to make good lower bounds on the albebraic and
arithmetic Hilbert functions in higher codimension and thus use the
Theorem of Minkowski effectively. 

Finally the argument, that the construed cycles with good approximation,
and bounded height of dimension $1$ defined over $\Z$ for different $D$'s
eventually differ is also borrowed from \cite{RW}. In the higher dimensionl
case however, one can not just use two one dimensional cycles to achieve the
result, but must prove that for certain one dimensional cycles $\alpha_D$, 
there also is one codimensional cycles, quasi defined over $\Z$ in a 
certain sense, that does not contain $\alpha_D$, and has also good
approximation properties. This is the principal part of the proof of
the second step mentioned above, and involves a rather complicated
combinatorics.

\section{The algebraic distance}

Let $\Pe^t_\Z$ be the projective space of dimensin $t$ over
$\spec \; \Z$, and effective cycles in $\Pe^t(\C)$ of codimensions $p$, 
and $q$ respectively. For $p+q \leq t+1$, and $X, Y$ properly intersecting or
in \cite{App1} the so called algebraic distance
\[ D(X,Y) \in \R, \]
is defined is known in the literature also as the 
height pairing of $X$ and $Y$ at infinity.
Further for $p+q \geq t+1$, and $X = \theta$ a point not contained int
the support of $Y$, there
are several essentially equivalent definitions of the algebraic distance
\[ D(\theta,Y) \in \R. \]
(\cite{App1}, section 4)

If further $|\cdot,\cdot|$ denotes the Fubini-Study ditance on $\Pe^t$ 
normalized in such a way that the maximal distance of two points is $1$,
there are the following Theorems for the algebraic distance.

\satz{Theorem I} \label{bezout}
Let $\CX,\CY$ be effective cycles 
intersecting properly, and $\theta$ a point in 
$\Pe^t(\C) \setminus (supp(X_\C\cup Y_\C))$. 
\begin{enumerate}

\item
There are effectively computable constants $c,c'$ only depending
on $t$ and the codimension of $X$ such that

\[ \deg(X) \log |\theta,X(\C)| \leq D(\theta,X) + c \deg X \leq 
   \log |\theta,X(\C)| + c' \deg X, \]
\item
If $f \in \Gamma(\Pe^t,O(D))_\Z$, and
$\CX = \di(f)$ is an effective cycle of codimension one,
\[ h(\CX) \leq \log |f_D|_{L^2} + D \sigma_t, \quad \mbox{and} \]
\[ D(\theta,X) + h(\CX) = \log |\la f| \theta \ra| + D \sigma_{t-1}, \]
where the $\sigma_i's$ are certain constants, and
$|\la f|  \theta \ra|$ is taken to be the norm of the evaluation of 
$f \in Sym^D(E) = \Gamma(\Pe^t,O(D))$ 
at a vector of length one representing $\theta$.

\item {\bf Metric B\'ezout Theorem}
For $p+q \leq t+1$,
assume that $\CX$, and $\CY$ have pure codimension $p$, and $q$ 
respectively.
There exists an effectively computable positive constant $d$,
only depending on $t$, and a map 
\[ f_{X,Y}: I \to \underline{\deg X} \times \underline{\deg Y} \]
from the unit interval $I$ to the set of natural numbers
less or equal $\deg X$ times the set of natural numers less or equal
$\deg Y$ such that $f_{X,Y}(0) = (0,0), f_{X,Y}(1) = (\deg X, \deg Y)$,
and the maps $pr_1 \circ f_{X,Y}:I \to \underline{\deg X}, 
pr_2\circ f_{X,Y}: I \to \underline{\deg Y}$
are monotonously increasing, and surjective, fullfilling: For every
$T \in I$, and $(\nu,\kappa) = f_{X,Y}(T)$, the inequality
\[ \nu \kappa \log |\theta,X+Y| + D(\theta,X.Y) + h(\CX . \CY) \leq \]
\[  \kappa D(\theta,X) + \nu D(\theta,Y) +
   \deg Y h(\CX) + \deg X h(\CY) + d \deg X \deg Y \]
holds. 
\item
In the situation of 3, if further
$|\theta,X+ Y| = |\theta,X|$, then
\[ D(\theta,X . Y) + h(\CX . \CY) \leq D(\theta,Y) +
   \deg Y h(\CX) + \deg X h(\CY) + d' \deg X \deg Y \]
with $d'$ a constant only depending on $t$.
\end{enumerate}
\end{Satz}

\proof
\cite{App1}, Theorem I.

Part 4 has an easy corollary:

\satz{Corollary}
\[ D(\theta,X . Y) + h(\CX . \CY) \leq \]
\[ max(D(\theta,X),D(\theta,Y)) +
   \deg Y h(\CX) + \deg X h(\CY) + d' \deg X \deg Y \]
\end{Satz}

Part three and four have a variant for arbitrary cycles over $\C$:
For $\CX,\CY$ defined over $\Z$, by \cite{App1}, Scholie 4.3,
\[ D(X,Y) = h(\CX.\CY) - \deg Y h(\CX) - \deg X h(\CY) + c \deg X \deg Y, \]
with $c$ a constant only depending on $p,q$, and $t$.
Replacing this into the formulas of the Theorem gives

\satz{Proposition} \label{bezext}
Let $Y$, and $Z$ be properly effective cycles over $\C$, and 
$\theta \in \Pe^t(\C) \setminus (supp Y \cup supp Z)$.
\begin{enumerate}

\item
There exists an effectively complutable constant $\bar{d}$, only
depending on $t$, and a map
\[ f_{Y,Z}: I \to \underline{\deg X} \times \underline{\deg Z}, \]
with properties as in the Theorem, such that for all $t \in I$, and
Then, with $\nu$, $\kappa = f_{Y,Z}(t)$,
\[ \nu \kappa \log |\theta,Y+Z| + D(\theta,Y.Z) + D(Y,Z) \leq \]
\[ \kappa D(\theta,Z) + \nu D(\theta,Y) + \bar{d} \deg Y \deg Z. \]

\item
If $|Z,\theta| \leq |Y, \theta|$, then
\[ D(\theta,Y.Z) + D(Y,Z) \leq D(\theta,Y) + \bar{d}' \deg Y \deg Z. \]

\item
\[ D(\theta,Y.Z) + D(Y,Z) \leq 
   \mbox{max}(D(\theta,Z),D(\theta,Y)) + \bar{d}' \deg Y \deg Z. \]
\end{enumerate}
\end{Satz}

\vspace{2mm}

Let now $\CY$ be defined over $\Z$, $f \in \Gamma(\Pe^t,O(D))_\Z$ with 
$f|_Y \neq 0$, define $f_Y^\bot$ as the orthogonal projection of $f$ modulo
the $D$-homogeneous part of the vanishing ideal $I_Y(D)$ of $Y$, and
$\CX := \di f, Z =Z_\C:= \di f_Y^\bot$. Then,
\begin{eqnarray*}
h(\CY .\di f) &=& D h(\CY) + \int_Y \log |f| \mu \\ \\
&=& D h(\CY) + \int_Y \log |f_Y^\bot| \mu \\ \\
&=& D h(\CY) + \deg Y \log \int_{\Pe^t} \log |f_Y^\bot| +
   D(Y,\di f_Y^\bot) \\ \\
&\leq&
   D h(\CY) + \deg Y \log |f_Y^\bot|_{L^2(\Pe^t)} + c D \deg Y, 
\end{eqnarray*}
by \cite{App1}, Proposition 6.3, with a constant $c$ only depending
on $t$. Further
\[ D(Y,Z) = \int_Y \log |f_Y^\bot| - \deg Y \int_{\Pe^t} \log |f_Y^\bot| =
   h(\CY . \di f) - D h(\CY) - \deg Y  \log |f|_{L^2(\Pe^t)}. \]
By \cite{App1}, 
\[ \int_{\Pe^t} \log |f_Y^\bot| \mu^t \geq \log |f_Y^\bot|_{L^2(\Pe^t)}
   - c_1 D, \]
and by the proof of \cite{App2}, Theorem 4.2, 
\[ \log |f_Y^\bot|_{L^2(\Pe^t)} \geq - \frac{h(\CY)}{\deg Y} D - c_2 D. \]
with $c_1,c_2$ only depending on $t$.
Hence, for every $\kappa$ with $0 \leq \kappa \leq \deg Y$,
\begin{eqnarray*}
\kappa D(Z,\theta) &=& \kappa\left( \log |\la f_Y^\bot|\theta \ra| -
   \int_{\Pe^t} \log |f_Y^\bot| \mu^t \right)  \\ \\
 &\leq& 
   \kappa \left( \log |\la f_Y^\bot|\theta \ra| - \log |f_Y^\bot|_{L^2(\Pe^t)}
   + c_1 D \right)  \\ \\
&\leq& \kappa \left( \log |\la f_Y^\bot|\theta \ra| + \frac{h(\CY)}{\deg Y} D +
   c_2 D \right)  \\ \\
&\leq& \kappa \log |\la f_Y^\bot|\theta\ra| + h(\CY) D + (c_1+c_2) D \deg Y. 
\end{eqnarray*}
Consequently, the previous Proposition implies

\satz{Propsosition} \label{bezexthe}
With the notations of Theorem \ref{bezout}, and $f, f_Y^\bot$ as above,
\begin{enumerate}

\item
\[ h(\CY. \di f) \leq D h(\CY) + \deg Y \log |f_Y^\bot| + \tilde{d} D \deg X.\]

\item
\[ \nu \kappa \log |\theta,Y+Z| + D(\theta,Y.Z) + h(\CY . \di f) \leq \]
\[ \kappa \log |\la f_Y^\bot|\theta\ra| + \nu D(\theta,Y) + 2 D h(\CY) + 
   \deg Y \log |f_Y^\bot|_{L^2(\Pe^t)} + \check{d} D \deg Y. \]

\item
If $|Z,\theta| \leq |Y, \theta|$, then
\[ D(\theta,Y.Z) \leq D(\theta,Y) + 2 D h(\CY) + 
   \deg Y \log |f_Y^\bot|_{L^2(\Pe^t)} + \bar{d}' \deg Y \deg Z. \]
If $|Z,\theta| \geq |Y, \theta|$, then
\[ D(\theta,Y.Z) \leq \log |\la f_Y^\bot|\theta\ra| + 2 D h(\CY) + 
   \deg Y \log |f_Y^\bot|_{L^2(\Pe^t)} + \bar{d}' \deg Y \deg Z. \]

\item
\[ D(\theta,Y.Z) \leq \]
\[ \mbox{max}(D(\theta,Z),\log |\la f_Y^\bot|\theta\ra|) + 
   2 D h(\CY) + \deg Y \log |f_Y^\bot|_{L^2(\Pe^t)}+\bar{d}' \deg Y \deg Z. \]
\end{enumerate}

\end{Satz}

One of the Definitions of algebraic distance is as follows:
$X \in Z^p_{eff}(\Pe^t_\C)$ be an effective of pure codimension one,
and $\Pe(W) \subset \Pe^t_\C$ a $p$- dimensional projective sub space
that intersects $X$ properly. If
\[ X . \Pe(W) = \sum_{x \in supp(X. \Pe(W))} n_x x, \]
where $n_x$ is the intersection multiplicity of $\Pe(W)$ and $X$ at $x$,
define $\rho_{\Pe(W}) := \sum_{x \in supp(X. \Pe(W))} n_x \log |x,\theta|$, and
the algebraic distance
\[ D_{pt} := \genfrac{}{}{0pt}{1}{\mbox{inf}}{\Pe(W)} \; \rho_{\Pe(W)}. \]
With these notations, and the ones of Theorem \ref{bezout}, one further has

\satz{Proposition} \label{bezmult}
If $\Pe(W)$ is a subspace of $\Pe^t$ where $\rho_{\Pe(W)}$ attains its
infinumum, and $M$ is the set of $y \in supp \Pe(W) \cap Y$ such that
$|Z,\theta| \leq |y,\theta|
(|X,\theta| \leq |y,\theta||\di f,\theta| \leq |y,\theta|)$
then
\[ D(Y . Z,\theta) + D(Y,Z) \leq 
   \sum_{y \in M} n_y \log |y,\theta| + d \deg X \deg Y. \]
\[ D(X . Y,\theta) + h(\CX . \CY) \leq 
   \sum_{y \in M} n_y \log |y,\theta| + \deg Y h(\CX) + \deg X h(\CY) + 
   d \deg X \deg Y. \]
\[ D(\di f . Y,\theta) \leq  
   \sum_{y \in M} n_y \log |y,\theta| + D h(\CY) +
   \deg Y \log |f_Y^\bot|_{L^2(\Pe^t)} + d \deg X \deg Y. \]
\end{Satz}

\proof
Reconsider the proof of Theorem \ref{bezout}.3 in \cite{App1}.
Let $\Pe(\bar{W})$ be a subspace that intersects $X$ properly, and
contains a point $z_0$ with $|z_0,\theta| = |Z,\theta|$. Then,
as in the proof of part 1. in the proof of \ref{bezout}.3 in \cite{App1},
\[ D_0((\theta,\theta), X \# Z) \leq
   \sum_{z \in Z \cap \Pe(\bar{W}), y \in Y \cap \Pe(W)} 
   n_z n_y \log |x \# z,(\theta,\theta)| + e_3 \deg Z \deg Y \leq \]
\[ \sum_{y \in M} n_{z_0} n_y
   \log |z_0 \# y,(\theta,\theta)| + e_3 \deg Z \deg Y , \]
which by \cite{App1}, Lemma 6.4 is less or equal
\[ \sum_{y \in M} n_y \log |y,\theta| + e_3 \deg Z \deg Y. \]
As $D(X,Z) \leq d \deg X \deg Z$ for some $d$,
the first inequality follows. 
From this the second inequality follows in the same way as par two of
the proof of Theorem \ref{bezout}.3 in \cite{App1}, again using
\[ D(\Pe(\Delta), X \# Y) = D(X,Y) + \bar{d} \deg X \deg Y = \]
\[ h(\CX . \CY) - \deg Y h(\CX) - \deg X h(\CY) -
   t \log 2 \deg X \deg Y  + \bar{d} \deg X \deg Y. \]
The second inequality follows in the same way as Proposition \ref{bezexthe}
above.

\satz{Definition}
Let $X$ be a quasiprojective  irreducible algebraic variety over a number 
field $k$, and $\sigma: k \to \C$ some imbedding. A point 
$\theta \in X_\sigma(k)$ is called generic if $\theta$ is contained in no 
proper algebraic subvariety $Y \subset X$ defined over $k$.
\end{Satz}

\vspace{2mm}

{\bf Remark:}
If $Y,Z$, cycles in $\Pe_\Z^t$ are defined over $\Q$, and 
$\theta\in \Pe^t(\C)$ is a generic point,
then automatically $\theta \notin supp Y \cup supp Z$, and the Propositions
of this section are applicable. However, Propositions
\ref{bezexthe}, and \ref{bezmult} will be have to be applied if possibly
$\la f_Y^\bot|\theta\ra = 0$. But it is easy to see that Theorem \ref{bezout}
and the other Propositions still hold if say 
$\theta \in supp \; X \setminus supp \;Y$ in Theorem \ref{bezout}:
If $\nu \kappa \neq 0$, then both sides in Theorem \ref{bezout}.3 are
$-\infty$, and thus the statement trivially holds. 
If $\nu \kappa = 0$, both sides of
the inequality are finite. To see that the inequality holds
let $(\theta_n)_{n \in \N}$ be a series of points in $\Pe(\C)$ such that
$\theta_n \notin supp \; X \cup supp \;Y$, 
$|\theta_n,X| < |\theta_n,Y|$ for all $n$, and
$lim_{n \to \infty} \theta_n = \theta$. Then, the inequality holds for each 
$\theta_n$ instead of $\theta$ and by continuity also for $\theta$.

\section{Hilbert functions}

A subscheme $\CX$ in $\Pe^t$ will be calle a subvariety if each
irreducible component has at least on $\bar{\Q}$-valued point
A global section 
$f \in \Gamma(\Pe^t_{\CO_k},O(D))$ is called primitive if it can not
be devided by any $a \in \CO_k$ which is not a unit. Global sections
will always be assumed to be primitive. The proper intersection
of a subvariety with the divisor of a primitive global section is
again a subvariety.

A subvariety $\CX \subset \Pe^t$ is called a locally complete intersection
if $\codim  \; \CX=r$, and there exist global sections $f_1,\ldots,f_r$, and
a Zariski open subset $U \subset \Pe^t$ such that
\[ \CX = \overline{V(f_1) \cap \cdots \cap V(f_r) \cap U}, \]
where $V(f_i)$ denotes the vanishing set of $f_i$.

Let $\CY \subset \CX \subset \Pe^t$ be algebraic subvarieties
with $\CY$ irreducible. $\CX$ is called a complete intersection at $\CY$,
if $\codim \CX =r$, and there are global sections $f_1,\ldots,f_r$ such 
that $\CX$ consists of the irreduciblecomponents of 
$V(f_1) \cap \cdots \cap v(f_r)$ that contain $\CY$.
The same notions can be defined for $k$-rational subvarieties of $\Pe^t$.
There is the easy Lemma:

\satz{Lemma} \label{lvdy}

\begin{enumerate}

\item
If $\CX \subset \Pe^t$ is a locally complete intersection, and 
$\CY \subset \CX$ is a subvariety of codimension zero, then $\CY$ is a
locally complete intersection.

\item
For any irreducible variety $\CY$; if $\CX$ is a locally complete intersection
at $\CY$, then $\CX$ is a locally complete intersection.

\item
If $\CX$ is a locally complete intersection at $\CY$, and $\CZ$ a subvariety
that contains $\CY$ and intersects $\CX$ properly, the union $\CW$
of the components of $\CX \cap \CZ$  that contain $\CY$ is a locally complete 
intersection at $\CY$.
\end{enumerate}
\end{Satz}

\satz{Proposition} \label{alHilbert}
Let $X$ be a subvariety of pure dimension $p$ in $\Pe^t_k$, and denote by 
\[ H_X(D) = \dim H^0(X,O(D)) \]
the algebraic Hilbert function.

\begin{enumerate}
\item
\[ H_X(D) \leq \deg X {D+p \choose p} \]

\item
If $X$ is locally complete intersection of
$p$ hyperplanes of degree $D_1, \ldots, D_p$, then
with $\bar{D}:= D_1+ \cdots + D_p -(t-p)$, and $D \geq \bar{D}$,
\[ H_X(D) \geq \deg X {D- \bar{D} + p \choose p}. \]
\end{enumerate}
\end{Satz}

\proof
1. \cite{Ch}, Theorem 1.

2. \cite{CP}, Corollaire 3.

An arithmetic bundle $\bar{M}$, for the purpose of this paper will just be a 
free $\Z$ module $M$ with a hermitian metric on $M_\C = M \otimes_\Z \C$.
The arthmetic degree $\wde (\bar{M})$ is then minus the logarithm of the
volume of the lattice $M$ in $M_\C$.

Let again $\Pe^t_\Z = \Pe(\Z^{t+1})$ be projective space of dimension $t$, and
\[ E_D := \Gamma(\Pe^t,O(D)). \]
As $E_D = Sym^D E_1$, which in turn equals the space of homogeneous
polynomiels of degree $D$ in $t+1$ variables, this lattice canonically 
carries the following metrics:

\begin{enumerate}

\item
The subspace metric $Sym^D E_1 \subset E_1^{\otimes D}$.

\item
The quotient metric $E_1^{\otimes D} \to Sym^D E_1$.

\item
The $L^2$-metric
\[ ||f||^2_{L^2(\Pe^t)} = \int_{\Pe^t(\C)} |f|^2 \mu^t, \]
where $\mu$ is the Fubini-Study metric on $\Pe^t$.

\end{enumerate}

Let $X$ be an effective cycle of pure dimension $s$ in $\Pe^t_\C$.
Then on 
\[ I_X(D) := \{ f \in H^0(\Pe^t,O(D)) | f|_X = 0 \}, \]
there are the restrictions of the norms $|\cdot|_{Sym}$, and
$|\cdot|_{L^2(\Pe^t)}$, and on
\[ F_X(D) = H^0(X,O(D)), \]
there is the quotiont norm $|\cdot|_{QN}$ induced by the canonical quotient map
\[ q_D: E_D \to F_D(X), \]
the $L^2(\Pe^t)$-norm
\[ | \cdot |_{L^2(\Pe^t)}: F_D(X) \to \R,  \quad
         \bar{f} \mapsto \inf_{q_D(f)=\bar{f}} \int_{\Pe^t} |f| \mu^t =
         \int_{\Pe^t} |g| \mu^t, \]
with $g \in I_X(D)^\bot$, and $q_D(g) = \bar{f}$.
and the $L^2(X)$-norm
\[ |\cdot|_{L^2(X)}: F_D(X) \to \R, \quad f \mapsto \int_X |f| \mu^p. \]

By convention $F_X(D)$ always denotes the $\Q$-or $\C$-vector space
of global sections, and $F_\CX(D)$ the corresponding lattice in $F_X(D)$.

\satz{Theorem} \label{arHilbert}
Let $\CX$ be an subvariety of dimension $s+1$ of $\Pe^t$, and denote by
\[ \hat{H}_\CX(D) := \wde (\bar{F}_\CX(D),|\cdot|_{L^2(\Pe^t)}), \quad
   \hat{\CH}_\CX(D) := \wde (\bar{F}_\CX(D),|\cdot|_{L^2(\CX)}) \]
the arithmetic Hilbert functions. 

\begin{enumerate}

\item
\[ \hat{\CH}_\CX(D)\leq \deg X
   \left( D h(\CX) + \frac12 (\log \deg X+2s \log D)\right) {D+s \choose s}. \]
\item
There is a positive constant $c_3$ only depending on $t$, 
such that
\[ \hat{H}_\CX(D) \leq \]
\[ \left( D h(\CX) + 
   \deg X D c_3 +
   \deg X (\frac12 \log \deg X + s \log D) \right) {D+s \choose D}. \]
Hence for $c_5 > c_3$, $\deg X$ at most a fixed
polynomial in $D$, and $D$ sufficiently large,
\[ \hat{H}_\CX(D) \leq D {D+s \choose D} (h(\CX)+c_5 \deg X). \]

\item
There is a constant $c_4$ only depending on $t$
\[ \hat{H}_\CX(D) \geq \wde (\overline{E_D \cap I_D(X)_\C^\bot}) \geq
   - \hat{H}_\CX(D)+ \frac{2 \sigma_t}{(t+1)!} D^{t+1} + 2 c_4 D^t \log D. \]
For $D$ sufficiently large, this is greater or equal
\[ - D {D+s \choose D} (h(\CX) + c_5\deg X). \]

\item
There are constants $c_1,c_2 > 0$ and $m \in \N$ only depending on $t$
such that if $\CX$ is an irreducible subvariety which is a locally complete 
intersection of $t-s$
hypersurfaces of degree $D_1, \ldots, D_{t-s}$, then with
$\bar{D} := D_1 + \cdots + D_{t-s}-(t-s)$, and $D \geq m \bar{D}$.
\[ \hat{H}_\CX(D) \geq (c_1 h(\CX) - c_2 \deg X) D^{s+1}. \]

\end{enumerate}
\end{Satz}

\proof
\cite{App2}, Theorem 4.1.





\section{Approximations in Projective Space}

\subsection{Approximation cycles}

An approximation cycle for a point $\theta \in \Pe^t_\C$ is an effective
cycle in $\Pe^t_\Z$ that has small algebraic distance to $\theta$ compared
with its degree and height. The precise definition, which with respect
to the ratio required between algebraic distance, height, and degree, 
is to a certain, but inessential degree arbitrary, will be given later. First
the general method to construct cycles with good approximation with respect
to $\theta$ will be exposed.

\subsubsection{Fundamental approximation techniques}

\satz{Theorem(Minkowski)} \label{Minkowski}
Let $\bar{M}$ be an arithmetic bundle over $\spec \; \Z$, and
$K \subset M_{\otimes \Z} \R$ any closed convex subset that is symmetric
with respect to the origin, and fullfills
\[ \log vol(K) \geq - \wde(\bar{M}) + \rk M \log 2. \]
Then $K \cap M$ contains a nonzero vector. In particular for
$c>0,\rk M \geq 2$, choosing a line $L \subset M_\R$ through the origin, and
taking $K$ as a rectangular paralellepiped with the edge parallel to $L$
of logarithmic length 
$c \; \wde \bar{M}$, and the edges orthogonal to $L$ with logarithmic 
length 
$-\frac{4c}{\rk M} \wde (\bar{M})+ 2\log 2$ centered at the origin, one sees
that there is a non zero lattice point $v \in M$ of logarithmic length
\[ \log |v| \leq - \frac{4c}{\rk M} \wde (\bar{M}) + 2 \log \rk M, \]
such that its projection to $L$ has logarithmic length at most 
$c \wde \bar{M}$.
\end{Satz}

\satz{Lemma}
Let $n \in \N$, and $v_1,\ldots,v_n \in \C^n$ be vectors such that
for every $i=1,\ldots,n$ the $i$th component of $v_i$ is nonzero. 
Then, there are numbers $m_1, \ldots, m_n \in \N$ with
$m_i \leq n$ such that each component of
\[ v= m_1 v_1 + \cdots m_n v_n \]
is nonzero. 
\end{Satz}

\proof
\cite{App2}, Lemma 2.6.

\vspace{2mm}

For any positive real number $a$, and an effective cylce $\CX$ in 
$\Pe^t$ define the $a$-size of $\CX$ as
\[ t_a(\CX) := a \deg X + h(\CX). \]

If $\CY$ is an irreducible subvariety of $\Pe^t$, identify
the arithmetic bundle 
$\bar{F}_\CY(D)=\overline{\Gamma(\Pe^t,O(D))/(I_\CY(D))}$, 
with the arithmetic bundle consisting of lattice
in $I_Y(D)^\bot$ consisting of all orthogonal projections of lattic 
points in $\Gamma(\Pe^t,O(D))$ ot $I_Y(D)^\bot$, and the induced
metric on $I_Y(D)^\bot$. If with this convention
\[ q_\CY: \Gamma(\Pe^t,O(D)) \to F_\CY(D) \]
denotes the quotient map, and $\CZ$ is an irreducibel subvariety of $\CY$,
then, for $f \in\Gamma(\Pe^t,O(D))$, 
\[ |q_Z(f)|_{L^2(\Pe^t)} \leq |q_Y(f)|_{L^2(\Pe^t)} \leq |f|_{L^2(\Pe^t)}. \]

\satz{Lemma} \label{klgs}
Let $\CX \subset \Pe^t$ be a locally complete intersection of global sections
$f_1,\ldots, f_s$ of degree $D_1 \leq \cdots \leq D_s$, $\theta \in \Pe^t(\C)$
a generic point, and define $\bar{D}:= D_1 + \cdots +D_s-s$. Then, with
$m,c_1,c_2$ the constants from Theorem \ref{arHilbert}.4, $D \geq \bar{D}$, 
$n \geq m$, $b \leq \mbox{min}(c_1,\frac1{2^{t+1}(t-s)!)}$, and
$a \geq c_2 2^{t+1}(t-s)!$,
there is an $f \in \Gamma(\Pe^t,O(nD))_\Z$, and an 
$\bar{f} \in \Gamma(\Pe^t,O(D))_\C$ such that for every irreducible
component $\CX_i$ of $\CX$ the quotions
$q_{X_i}(f),q_{X_i}(\bar{f})$ are nonzero,
and for every subvariety $\CZ \subset \Pe^t$ contained in each irreducible 
component of $\CX$, the quotients
$q_Z(f)$ and $q_Z(\bar{f})$ coincide and further
\[ \log |f_Z^\bot| = \log |\bar{f}_Z^\bot| \leq 6anD, \quad \mbox{and} \quad
   \log |\la \bar{f} |\theta\ra| \leq - b t_a(\CX_{min}) (nD)^{t+1-s}, \]
where $\CX_{min}$ denotes an irreducible component of $\CX$ with
minimal $a$-size.
\end{Satz}

\proof
Assume first $\CX = \CY$ is irreducible. 
By Proposition \ref{alHilbert}.,
\[ (nD)^t \geq \deg Y (nD)^{t-s} \geq
 \deg Y {nD+t-s \choose t-s} \geq H_Y(nD) = \rk F_Y(nD) \geq \]
\begin{equation} 
\deg Y {nD - \bar{D} + t-s \choose t-s} \geq\frac{\deg Y}{2^{t-s}(t-s)!}
 (nD)^{t-s},
\end{equation}
since $m \geq 2$.
and by Proposition \ref{arHilbert}.4,
\begin{equation}
- \wde(\bar{F}_\CY(nD)) = - \hat{H}_{\CY}(nD)
\leq (-c_1 h(\CY) + c_2 \deg Y) (nD)^{t+1-s}.
\end{equation}
Let next $L_\theta \subset F_Y(nD)_\C$ be the one dimensional
subspace orthogonal to the kernel $ker_\theta$ of the evaluation map
from $F_Y(nD)_\C$ to the stalk of $O(nD)$ at $\theta$, and
$K$ the rectangular parallelepiped with logarithmic
length of edge parallel to $L_\theta$ equal to 
\[ -b h(\CY) (nD)^{t+1-s}- b a \deg Y (nD)^{t+1-s}
   = - b t_a(\CY) (nD)^{t+1-s}, \]
and  logarithmic length of edges parallel to $ker_\theta$ equal to
$4anD$.
By the choice of $a$, and $b$, and the estimates on the algebraic and
arithmetic Hilbert functions above, we have
\begin{eqnarray*}
\log vol(K) &=& 4 a n D (\rk F_Y(nD)-1) 
   - b h(\CY) (nD)^{t+1-s}-b a\deg Y(nD)^{t+1-s} \\ \\
&\geq& \frac{4a}{2^{t-s+1}(t-s)!} \deg Y (nD)^{t+1-s}
   - b h(\CY) (nD)^{t+1-s} \\
&& \hspace{55mm} -b a\deg Y(nD)^{t+1-s} \\ \\
& \geq & \left( \left( \frac{4a}{2^{t+1-s}(t-s)!} - ba \right) \deg Y
           - b h(\CY) \right) (nD)^{t+1-s} \\ \\
& \geq& (c_2 \deg Y - c_1 h(\CY)) (nD)^{t+1-s} + \deg Y (nD)^{t-s} \log 2\\ \\
&\geq &- \wde(\bar{F}_\CY(nD)) + \rk F_Y(nD) \log 2. 
\end{eqnarray*}
Hence, by the Theorem of Minkowski, $K$ contains a non zero lattice point
$\bar{f}$, that is
\[ \log |\bar{f}| \leq 4anD + \frac12 \log \rk F_Y(nD) \leq
   4anD + \frac12 \log (nD)^t \leq 5anD \]
for $D$ sufficently big, and
\[ \log |\la \bar{f}|\theta\ra| \leq -b h(\CY) (nD)^{t+1-s}- 
   b a \deg Y (nD)^{t+1-s} = - b t_a(\CY) (nD)^{t+1-s}. \]

Choose $f \in \Gamma(\Pe^t,O(nD))_\Z$ as a representative of $\bar{f}$.
Then for $\CZ \subset \CY$, clearly $q_Z(f) = q_Z(\bar{f})$, and
\[ \log |f_Z^\bot| \leq |f_Y^\bot| = |\bar{f}| \leq 5anD, \]
proving the Lemma in case $\CX$ is irreducible.


Assume now $\CX$ is not irreducible, and let
$\CX = \CX_1 \cup \cdots \cup \CX_l$ be the decomposition of $\CX$ into
irreducible components, and $\CZ$ a variety contained in each
$\CX_i, i=1,\ldots,l$. We have $l \leq \deg Y \leq (nD)^s$.
By the previous argument for each $i$ there is
a nonzero $\bar{f}_i \in F_{\CX_i}(nD)$, and an
$f_i \in \Gamma(\Pe^t,O(D))_\Z$ such that $q_{X_i}(f_i) = \bar{f}_i$,
\[ \log |(f_i)_Z^\bot| \leq \log |(f_i)_{X_i}|^\bot \leq 5 a n D, \]
and
\[ \log |\la \bar{f}_i|\theta \ra| \leq 
   - b t_a(\CX_i) (nD)^{t+1-s} \leq
   - b t_a(\CX_{min}) (nD)^{t+1-s}. \]
By the previous Lemma there are natural numbers $k_i, i=1, \ldots l$
with $\log k_i \leq \log l \leq s \log D$ such that
\[ \bar{f} = \sum_{i=1}^l k_i \bar{f}_i \]
is nonzero on each $X_i$, and the same holds for
\[ f = \sum_{i=1}^l k_i f_i. \] 
Because of $q_{X_i}(f_i) = \bar{f}_i$, and $Z \subset X_i$ for every
$i=1, \ldots, l$, we have $q_Z(f) = q_Z(\bar{f}$, hence
\[ \log |f_Z^\bot| = \log |\bar{f}_Z^\bot| \leq \log |\bar{f}| \leq
   \mbox{max}_{i=1,\ldots,l} \log |\bar{f}_i| + s \log D \leq
   5 a n D + s \log D \leq 6anD, \]
for $D$ sufficiently large, and
\[ \log |\la \bar{f}|\theta \ra| \leq 
   - b t_a(\CX_{min}) (nD)^{t+1-s} + s \log D, \]
proving the Lemma.

\vspace{2mm}

{\bf Remark:}
It is easily seen that the Lemma holds with $\bar{b} = \frac nm b$ instead
of $b$ if $n/m$ is still bigger than $m$. This fact will be assumed when 
the Lemma is aplied.

\vspace{2mm}

A strategy for constructing effective cycles of arbitrary codimension
of bounded height and degree with small algebraic distance 
to a given generic point $\theta$, that is approximation cycles, could now
be the following. Use the previous Lemma to find a vector
$f_1 \in \Gamma(\Pe^t,O(D))$ of bounded length such that
$|\la f_1|\theta\ra|$ is small, define $\CX_1:= \di f$, and use 
Theorem \ref{bezout}.2 to derive that $\CX_1$ has small algebraic distance
with respect to $\theta$. It is then easy to prove that some irreducible 
component $\CY_1$ of $\CX_1$ is also an approximation cycle. 
Then, there are two possibilities to go on. 

The first one uses the previous Lemma to find
an $f_2$ of bounded length and degree that is nonzero on $\CY_1$ and has
small $|\la f_2|\theta\ra|$, and the metric B\'ezout Theorem
to prove that $\CX_1:= \CY_1 . \di f$ has small algebraic distance
with respect to $\theta$, and again some irreducible cycle $\CY_2$ of
$\di f . \CY_1$ will have good approximation with respect to $\theta$. 
This possiblity is chosen in \cite{Ph}, only
with the weaker estimate in \ref{arHilbert}.3 for the arithmetic Hilbert 
function
thereby supplying a weaker approximation than possible with the previous 
Lemma which rests on the estimate \ref{arHilbert}.4.
The problem with this approach is that there is no guarantee that
the successive intersections $\CY_3,\CY_4,\ldots$ (i.\@ e.\@ from 
codimension 3 onwards) are locally complete 
interesections of bounded degree, and hence the previous Lemma is not
applicable.

The second possiblity is to use the Lemma for $\CX_1$ to obtain an 
$f_2$ intersecting each irreducible component of $\CX_1$ properly and having
small $|\la f_2|\theta \ra|$; one does then not leave the realm of
locally complete intersections. The problem with this approach is that if
the irreducible component with minimal $a$-size 
$\CX_{min}$ of $\CX_1$ has very small $a$-size,
the estimate on $\log |\la f_2|\theta \ra|$ one obtains is not very good.
If $\CX_{min}$ itself has small algebraic distance to $\theta$, this 
doesn't matter, because the metric B\'ezout Theorem will still give a
good estimate, but if $\CX_{min}$ has big albebraic distance one is in 
a trap. To get out, one has to prove that one does not come along a
$\CX_{min}$ with
big algebraic distance if one applies a somewhat refined procedure, which
involves not only constructing approximation cycles of higher codimension
from ones of lower codimension in the way just sketched but also, if it 
should not be possible to construct an approximation cycle of higher 
codimension, to obtain one which has
lower codimension but better approximation properties  than any of the
ones so far constructed. To prove that this road finally leads to
cycles of codimension $t$ involves rather complex combinatorics which will be
presented in the next subsection, but first we prove how to construct the
approximation cycle of lower codimension with better approximation.

\satz{Lemma} \label{hoch}
Let $r<s<t$ be natural numbers, 
$\CY,\CX = \CX_1 \cup \cdots \cup \CX_l$ varieties in $\Pe^t$ of codimensions 
$s,r$,
with $\CY$ and every component $\CX_i , i =1,\ldots,l$ of $\CX$ containing 
$\CY$. Let further $\CX_{min}$ be the an irreducible component of $\CX$ with
minimal $a$-size, $n$ a natural number with $n << D$, and
$b_s, b_{s+1},\bar{b}_s \leq \bar{b}_r$ positive real numbers
with $b_{s+1} \leq b_s/c$. 
Finally $f \in \Gamma(\Pe^t,O(nD))_\Z$,
$\bar{f} \in \Gamma(\Pe^t,O(nD))_\C$ such that
$\bar{f}_Y^\bot = f_Y^\bot \neq 0$, and 
$\bar{f} = \sum_{i=1}^l \bar{f}_i$, $f = \sum_{i=1}^l f_i$ such that  
$\bar{f}$, and $f$ are nonzero
on each irreducible component of $\CX$, and
$f_i,\bar{f}_i$ are nonzero on $\CX_i$.
If these date fullfill the inequalities
\[ \log |f_Y^\bot| = \log |\bar{f}_Y^\bot| \leq 6anD, \quad
   \log |(f_i)_{X_i}^\bot| \leq 6anD, \quad
   \log |\la \bar{f}_i|\theta\ra| \leq - \bar{b}_r t_a(\CX_i) D^{t+1-r}, \]
\[ D(Y,\theta) \leq - b_s t_a(\CY) D^{t+1-s}, \quad
   b_s t_a(\CY) \geq 2 \bar{b}_r t_a(\CX_{min}) D^{s-r}, \]
then $t_a(\di f . \CY) \leq 7 n t_a(\CY)$, and either 
$D(\di f .Y,\theta) \leq -b_{s+1} t_a(\di f.\CY) D^{t-s}$, or 
\[ D(\CX_{min}, \theta) \leq - \bar{b}_r t_a(\CX_{min}) D^{t+1-r}. \]
\end{Satz}

\proof
The claim on $t_a(\di f . \CY)$ follows from the algebraic B\'ezout
Theorem, and Proposition \ref{bezexthe}.
For the other claim, make the case distinction

{\sc Case 1: 
$\deg Y \geq \frac{b_s t_a(\CY)}{2\bar{b}_r t_a (\CX_{min})D^{s-r}}$.}

Choose $t$ in Proposition \ref{bezexthe}.2 such that 
$\mu = \frac{b_s t_a(\CY)}{2\bar{b}_r t_a (\CX_{min})D^{s-r}} \leq \deg Y$.
Then, by the assumptions on $D(Y,\theta)$, and $|\la\bar{f_i}|\theta\ra|$,
and the fact $h(\di f . \CY) \geq 0$,
\[ \mu \nu \log |Y,\theta| + D(\di f. Y, \theta) \leq \]
\[ \mu \log |\la \bar{f}|\theta\ra| + \nu D(Y,\theta) + 6 n D h(\CY) + 
   \deg Y 2 a n D + d n D \deg Y \leq \]
\[ - \mu \bar{b}_r t_a(\CX_{min}) D^{t+1-r} - \nu b_s t_a(\CY) D^{t+1-s} +
   7 n t_a(\CY) D = \]
\[ - \frac12 b_s t_a(\CY) D^{t+1-s}- \nu b_s t_a(\CY) D^{t+1-s} +
    7 n t_a(\CY)D = \]
\[ - \frac12 b_s t_a(\CY) D^{t+1-r}-\mu \nu 2\bar{b}_r t_a(\CX_{min}) 
   D^{t+1-r} + 7 n t_a(\CY)D. \]
Since $s \leq t-1$, and $n<<D$, either
\[ D(\di f . Y,\theta) \leq
   -\frac14 b_s t_a(\CY)D^{t+1-s}-\frac{\nu}2 b_s t_a(\CY) D^{t+1-s}
   \leq  -\frac14 b_s t_a(\CY)D^{t+1-s} \leq \]
\[ -b_{s+1} t_a(\di f.\CY) D^{t-s}, \]
and the first possibility holds, or 
\[ \log |Y,\theta| \leq  - \bar{b}_r t_a(\CX_{min}) D^{t+1-r} - 
   \frac1{4\mu\nu} b_s t_a(\CY)D^{t+1-s} \leq 
   - \bar{b}_r t_a(\CX_{min}) D^{t+1-r}, \]
which by Theorem \ref{bezout}.1, implies, because $Y$ is contained in $X_{min}$
\[ D(X_i,\theta) \leq \log |X_{min},\theta|+ c_1 \deg X_{min} 
   \leq \log |Y,\theta| + c_1 \deg X_{min}\leq \] 
\[ - \bar{b}_r t_a(\CX_{min}) D^{t+1-r} + O(t_a(\CX_{min})). \]
and hence $\CX$ the second possibility holds.

{\sc Case 2:
$\deg Y \leq \frac{b_s t_a(\CY)}{2\bar{b}_r t_a (\CX_{min})D^{s-r}}$.}
In this case, by Theorem \ref{bezout}.1, and the assumption on
$D(Y,\theta)$,
\[ \log |Y,\theta| \leq \frac1{\deg Y} D(Y,\theta) \leq
   - \frac {b_s}{\deg Y} t_a(\CY) D^{t+1-s} \leq
   - 2 \bar{b}_r t_a(\CX_{min}) D^{t+1-r}, \]
and the second possibility follows in the same way as above.

\subsubsection{The combinatorics}

Let $c_1,c_2,m$ be the constants from Theorem \ref{arHilbert}.4,
$d$ the constant from \ref{bezexthe}.1,
fix a real number $a>>0$ and a number $N \in \N$, and define the
constants
\begin{equation} \label{const}
\begin{array}{llll}
  \bar{n}_1 := 1, &\quad& 
   \bar{n}_s := Nm (2(1+m))^{s-2}, \quad &2 \leq s \leq t \\ \\
 n_1 := \bar{n}_1, &\quad&
   n_{s+1} := \bar{n}_{s+1} n_s, \quad &1 \leq s \leq t-1, \\ \\
a_1 := a, &\quad& a_{s+1} := 
        6 a \bar{n}_{s+1} n_s + a_s \bar{n}_{s+1} + d n_s \bar{n}_{s+1}, 
                       \quad &1 \leq s \leq t-1, \\ \\
&&\bar{m}_s := \bar{n}_{s+1} + 8 \bar{n}_s +1, \quad& 1 \leq s \leq t, \\ \\
 && m_s := \prod_{i=s}^t \bar{m}_i, \quad &1 \leq s \leq t, \\ \\
 \bar{b}_1:= \frac{N^t}{(t!)}, &\quad& \bar{b}_s:= N^t
     \mbox{min} \left(\frac1{2^{t+1}(t-s)!},c_1 \right), \quad 
                       &2 \leq s \leq t, \\ \\
  b_1 := \bar{b}_1, &\quad&
   b_s := \mbox{min}
    \left(\frac{b_{s-1}}{16 \bar{n}_s},\mbox{min}_{r=1,\ldots,s-1}
    \bar{b}_{r} \frac{m_s a_r}{m_r a_s} \right), \quad &2 \leq s\leq t. 
\end{array}
\end{equation}
Observe that all constants $n_i,\bar{n}_i,\bar{m}_i, m_i b_i, \bar{b}_i$, and
all quotients $a/a_i, a_i/a, i=1, \ldots,t$ are bounded
by a constant only depending on $t$ and $N$.

These constants are used to make a number of Definitions that also
illustrate their role. For the purpose this paper,
only the case $N=1$ is important, but for later applications the 
general case will be needed. For this reason it may be advisable on first
reading to restrict to the case $N=1$.

\satz{Definition and Lemma} \label{sucsch}
A chain of irreducible subvarieties
\[ \Pe^t = \CY_0 \supset \CY_1 \supset \cdots \supset \CY_s \]
is called a successive $D$-intersection (relative to $N$) if for
each $i=1, \ldots, s$ the corresponding variety $\CY_i$ is an irreducible 
component
of the proper intersection of $\CY_{i-1}$ with a global section
$f\in\Gamma(\Pe^t,O(\bar{n}_i D))$, and has thereby codimension one.
For a successive intersection the inequality
$\deg Y_i \leq n_i D^i$ holds. The chain is called a successive
$(D,a)$-intersection if additionally $h(\CY_s) \leq a_s D^s$.
\end{Satz}

\proof
Clearly, $\deg \CY_0=\deg \Pe^t = 1 = D^0$.
Further if $\CY_{i+1}$ is an irreducible component of $\di f . \CY_i$,
with $f \in \Gamma(\Pe^t,O(D))$,
then 
\[ \deg Y_{i+1}\leq \deg f \deg Y_i \leq \bar{n}_{i+1} D n_i D^i = 
   n_{i+1} D^{i+1}. \]


\satz{Definition and Lemma} \label{lvdd}
A subvariety $\CX \subset \Pe^t$ of constant codimension $s$ is called
a locally complete $D$-intersection (relative to $N$), 
if $\CX$ is a locally complete
intersection of $f_1, \ldots, f_s$ with $\deg f_i \leq \bar{n}_i D$.
$\CX$ is called a locally complete $D$- intersection at $\CY$ if 
$\CX$ is a locally complete intersection of $f_1, \ldots,f_s$ at $\CY$, with
$\deg f_i \leq \bar{n}_iD$. 
Further $\CX$ is called a locally complete $(D,a)$-intersection
(at $\CY$) if additionally $h(\CX) \leq a_s D^s$.
A locally complete $(D,a)$-intersection $\CX$ (at $\CY$) fullfills
$\deg \CX \leq n_s D^s$.
\end{Satz}

\satz{Definition and Lemma} \label{dh}
For an effective cycle $\CX$ of pure codimension $r$ the number
\[ t_a(\CX) m_r D^{\dim \CX} \]
will be called the dimensional $a$-size of order $D$. 
If $f \in \Gamma(\Pe^t,O(\bar{n}_{r+1} D))$ intersects a
subvariety $\CY$ of $\CX$ that has pure codimension $s > r$ in $\Pe^t$ 
properly,and fullfills $\log |f_{Y}^\bot| \leq 6 a \bar{n}_{r+1} D$, and 
$\CZ$ is any irreducible component of $\CY . \di f$, then
\[ t_a(\CZ) m_{r+1} D^{\dim \CZ} \leq t_a(\CX_i . \di f) m_{s+1}
   D^{\dim (\CX . \di f)} < t_a(\CY) m_s D^{\dim \CY}. \]
\end{Satz}

\proof
The first inequality is obvious. Further, by the algebraic B\'ezout
Theorem, and Proposition \ref{bezexthe}.1, with $a >d$,
\begin{eqnarray*}
t_a(\CZ) &\leq&
   t_a(\CY. \di f)  = a \deg (\CY . \di f) + h(\CY . \di f) \\ 
&\leq & a \bar{n}_{r+1} D \deg Y + 6 a \bar{n}_{r+1} D \deg Y + 
   \bar{n}_{r+1} D h(Y) + d \bar{n}_{r+1} D \deg Y \\ 
&\leq & a \bar{n}_{r+1} D h(\CY) + 7 a \bar{n}_{r+1}D \deg Y \\ 
&< & a D \bar{m}_r  \deg Y + \bar{m}_r D h(\CY) \\ 
&=& t_a(\CY) \bar{m}_r D \leq 
   t_a(\CY) \bar{m}_s D, 
\end{eqnarray*}
since $\bar{m}_r < \bar{m}_s$.
Consequently,
\[ t_a(\CZ) m_{s+1} D^{\dim (\CY . \di f)} =
   t_a(\CZ) \frac{m_s}{\bar{m}_s} D^{\dim (\CY . \di f)} \leq \]
\[ t_a(\CY . \di f) \frac{m_s}{\bar{m}_s} D^{(\dim \CX .\di f)}<
   t_a(\CY) m_s D^{\dim \CY}. \]
 
\vspace{2mm}

The dimensional size of a cycle $\CY$ measures the approximation power
of the space of global sections on $\CY$ with bounded length and degree.

\satz{Definition and Lemma} \label{appcyc}
For a given $\theta \in \Pe^t(\C), a>0$, and an effective cycle $\CX$, 
of pure codimension $s$ in $\Pe^t$, such
that the support of $\CX$ does not contain $\theta$, define the weighted
algebraic distance of $\theta$ to $\CX$ 
\[ \varphi_{a,\theta}(\CX) := \frac{D(\theta,X)}{t_a(\CX)}. \]
The cycle $\CX$ is called an aproximation cycle, or $a$-approximation 
cycle of order $D$ relative to $N$ for $\theta$, iff
\[ \deg X \leq n_s D^s,\quad h(\CX)\leq a_s D^s, \quad \mbox{and} \quad
   \varphi_{a,\theta} (\CX) \leq -b_s D^{t+1-s}, \]
If $\CX_l, l= 1, \ldots, L$ are $L$ effective cycles, then
\[ \genfrac{}{}{0pt}{1}{\mbox{min}}{l=1,\ldots,L} 
   (\varphi_{a,\theta}(\CX_l)) \leq
   \frac1L \varphi_{a,\theta} \left(\sum_{l=1}^L \CX_l \right). \]
Hence, if a sum of $L$ effective cycles is an approximation cycles, 
at least one of the summands is likewise.
In particular
\[ \varphi_{a,\theta}(n \CX) = \varphi_{a,\theta}(\CX), \]
for any natural number $n$, and any effective cycle $\CX$, which
in turn implies that cycles with different weighted algebraic distance
can not be multiples of the same irreducible variety.
If $\CX$ is an approximation cycle for $\theta$, then
\[ \log |\theta,X| \leq - b_s D^{t+1-s} +c, \]
where $c$ is the constant from \ref{bezout}.1
\end{Satz}

\proof
The first claim follows from the additivity, of the degree, height, and
algebraic distance, and from the nonpositivity of the algebraic distance, 
and the nonnegativity of the degree, and height.
For the second claim, since $h(\CX) \geq 0$, Theorem \ref{bezout}.1 implies
\[ \log |\theta,X| \leq \frac{D(\theta,X)}{\deg X} + c \leq
   \frac{D(\theta,X)}{a \deg X + h(\CX)} + c = \]
\[ \varphi_{a,\theta}(\CX) + c \leq - b_s D^{t+1-s} +c .  \]

\vspace{2mm}

Denote now by $\bar{C}_{t,D}$ be the set of $(D,a)$-approximation chains, 
that is the set of all sucessive $(D,a)$-intersections
\[ \CC: \quad \Pe^t = \CY_0 \supset \cdots \supset \CY_s \]
such that $\CY_s$ is a $(D,a)$-approximation cycle. 
Denote further by $C_{t,D}$ the set of equivalence classes
in $\bar{D}_{t,D}$ where
two chains are said to be equivalent iff they have the same end term. Of
course, $C_{t,D}$ is then just the set of approximation cycles that appear as
succesive intersection. However, some of the following proofs will have a more
lucid appearance
if one views the set $C_{t,D}$ as a set of equivalence classes of chains
of cycles rather than a set of cycles. 
On the other hand it will not be necessary to distinguish between an
element of $C_{D,t}$ and one of its represantatives in $\bar{C}_{D,t}$, and
therefore $\bar{C}_{D,t}$ will not appear henceforth.

On $C_{t,D}$ furhter define the following relation. For two chains
\[ \CC: \quad \Pe^t = \CY_0 \supset \cdots \supset \CY_s, \quad \quad
   \bar{\CC}: \quad \Pe^t = \CY_0 \supset \cdots \supset \CY_{\bar{s}}, \]
the relation $\CC \prec \bar{\CC}$ holds iff there is a sequence of
approximation chains
\[ \CC=\CC_0, \CC_1, \ldots, \CC_{k-1}, \CC_k = \bar{\CC}, \]
such that for each pair $\CC_l,\CC_{l+1}, l=0, \ldots, k-1$ 
\[ t_a(\CY) m_{t+1-\dim \CY}  D^{\dim \CY} < 
   t_a(\CX) m_{t+1-\dim\CX} D^{\dim \CX}, \]
where $\CY,\CX$ are the end terms of $\CC_l$, and $\CC_{l+1}$ respectively,
and one of the following conditions holds

\begin{enumerate}

\item
$\CC_{l+1}, \CC_l$ look like
\[ \CC_{l+1}: \quad \CY_0 \supset \cdots \supset \CY, \quad \quad
   \CC_l: \quad \CY_0 \supset \cdots \supset \CY \supset \CX. \]

\item
$\CC_{l+1}, \CC_l$ lool like
\[ \CC_{l+1}: \quad \CY_0 \supset \cdots \supset \CY, \quad \quad
   \CC_l: \quad \CY_0 \supset \cdots \supset \CX, \]
with $\CX$ containing $\CY$.
\end{enumerate}

This relation is obviously transitive, and also antireflexive, since for chains
$\CC, \bar{\CC}$ with end terms $\CY,\bar{\CY}$ the relation
$\CC \prec \bar{\CC}$ implies
\[ t_a(\CY) m_{t+1-\dim \CY}D^{\dim \CY} < t_a(\bar{\CY}) 
   m_{t+1- \dim \bar{\CY}} D^{\dim \bar{\CY}}, \]
and hence $\CY \neq \bar{\CY}$ implying $\CC_1 \neq \CC_2$. 

Since for each $D \in \N$, there are only finitely many 
subvarieties $\CZ$ with $\deg Z \leq n_{\dim Z} D^{\dim Z}$, 
and $h(\CZ) \leq a_{\dim Z} D^{\dim Z}$, the set $C_{D,t}$ is finite,
and hence there are minimal approximation chains with respect to the 
relation $\prec$; i.\@ e.\@ there is at least one $\bar{\CC} \in C_{t,D}$
such that there is no $\CC \in C_{t,D}$ with $\CC \prec \bar{\CC}$.



\subsubsection{Existence of Approximation cycles}

The fundamental Theorem on approximation cycles is the following:

\satz{Theorem} \label{appcycle}
For each sufficiently big $a$, and very $D >> 0$ the minimal
chains in $C_{t,D}$ with respect to the relation $\prec$
have length $t$. In particular, there is a 
successive $(D,a)$ intersection $\CY_D$ of pure codimension $t$
which is an approximation cycle of order $D$.
\end{Satz}

\satz{Lemma} \label{lvdrsch}
Let $\CY \subset \Pe^t$ be a succesive $(D,a)$-intersection 
of codimension $s \leq t-1$, and
\[ \CX = \overline{f_1 \cap \cdots \cap f_r \cap U} \] 
a locally complete $D$-intersection at $\CY$ of codimension $r < s$. 
Then, there are
global sections $g \in \Gamma(\Pe^t,O(\bar{n}_{r+1} D))_\Z$
$\bar{g} \in \Gamma(\Pe^t,O(\bar{n}_{r+1} D))$, such that
$g_Y^\bot = \bar{g}_Y^\bot \neq 0$, the restrictions of $f,\bar{f}$ to
every irreducible component of $\CX$ are nonzero, and 
\[ \log |g_Y^\bot| \leq 6 a \bar{n}_{r+1} D, \quad
   \log |\la \bar{g}|\theta \ra| \leq - 
   \bar{b}_{r+1} t_a(\CX_{min}) D^{\dim \CX}, \]
and the restrictions of $f$, and $\bar{f}$ to every irreducible
component of $X$ are non zero.

\end{Satz}

\proof
Since $\CX$ is a locally complete $(D,a)$-intersection of the 
$f_1, \ldots, f_r$, we have
\[ \bar{D} := \sum_{i=1}^r \deg f_i- r \leq D \sum_{i=1}^r \bar{n}_i =
   \leq \frac{\bar{n}_{r+1}}m D. \]
Since, by definition of a locally complete intersection at $\CY$ the
variety $\CY$ is contained in every irreducble component of $\CX$, the 
claim follows from Lemma
\ref{klgs} with $n = \bar{n}_{r+1}, b = \bar{b}_{r+1}$.

\satz{Corollary} \label{cor}
Let $\CY$ be a successive $(D,a)$-intersection of
codimension $s \leq t-1$. Then, there is $\CX$ a 
locally complete $D$-intersection at $\CY$ of codimension $r \leq s$, and 
global sections $g \in \Gamma(\Pe^t,O(\bar{n}_{r+1} D))_\Z$,
$\bar{g} \in \Gamma(\Pe^t,O(\bar{n}_{r+1} D))$, such that
$g_Y^\bot = \bar{g}_Y^\bot\neq 0$, and 
\[ \log |g_Y^\bot| \leq 6 a \bar{n}_{r+1} D, \quad
   \log |\la \bar{g}|\theta \ra| \leq - 
   \bar{b}_{r+1} t_a(\CX_{min}) D^{\dim \CX}, \]
and the restrictions of $f$, and $\bar{f}$ to every irreducible
component of $X$ are non zero.
\end{Satz}

\proof
We proof by complete induction that $\CY$ is contained in a 
of locally complete $D$-intersection $\CX$ at $\CY$ 
of codimension $r$ at $\CY$,
or there is a global section of the form specified in the Corollary.
Of course this claim entails the Corollary.
Clearly the claim holds for $r=0$ with $\CX = \Pe^t$.
So assume there is a locally complete $D$-intersection $\CX_r$ at $\CY$
of codimension $r$.
By the previous Lemma, there are global sections
$g \in \Gamma(\Pe^t,O(\bar{n}_{r+1} D))_\Z$
$\bar{g} \in \Gamma(\Pe^t,O(\bar{n}_{r+1} D))$, such that
$g_Y^\bot = \bar{g}_Y^\bot$, and 
\[ \log |g_Y^\bot| \leq 6 a \bar{n}_{r+1} D, \quad
   \log |\la \bar{g}|\theta \ra| \leq - 
   \bar{b}_{r+1} t_a(\CX_{min}) D^{\dim \CX}, \]
and the restrictions of $f$, and $\bar{f}$ to every irreducible
component of $X$ are non zero.
If $f$ has nonzero restriction to $Y$ the Corollary follows. 
If the restriction of $f$ to $Y$ is zero, 
since $\deg g = \bar{n}_{r+1}$, the union of the irreducible components
of $\CX . \di f$ that contain $\CY$ by Lemma \ref{lvdy}, is a locally complete 
$D$-intersection of codimension $r+1$ at $\CY$.

\satz{Lemma} \label{zwischen}
For every irreducible $(D,a)$-approximation cycle $\CY$ of codimension
$s \leq t-1$, belonging to an approximation chain $\CC$, either there is a 
$(D,a)$-approxi- mation chain
\[ \CC_1: \quad \Pe^t = \CY_0 \supset \cdots \supset \CX, \]
such that $\CX$ is a locally complete $D$-intersection at $\CY$, and
$\CC_1 \prec \CC$, or there is a $(D,a)$-intersection chain
\[ \CC_2: \quad \CY_0 \supset \cdots \supset \CY \supset \CZ \]
with $\CC_2 \prec \CC_1$.
\end{Satz}

\proof
Let $\CX = \CX_1 \cup \cdots \cup \CX_l$ be the locally complete
$D$-intersection at $\CY$, and
$g \in \Gamma(\Pe^t,O(\bar{n}_{r+1} D))_\Z$,
$\bar{g} \in \Gamma(\Pe^t,O(\bar{n}_{r+1} D))_\C$, the global section
such that $g_Y^\bot = \bar{g}_Y^\bot\neq 0$, and 
\[ \log |g_Y^\bot| \leq 6 a \bar{n}_{r+1} D, \quad
   \log |\la \bar{g}|\theta \ra| \leq - 
   \bar{b}_{r+1} t_a(\CX_{min}) D^{\dim \CX} \]
from the Corollary.

By the Theorem of B\'ezout,
\[ \deg (Y . \di f) = \deg Y \bar{n}_{r+1} D \leq
   n_s \bar{n}_{r+1} D^{s+1} \leq n_s n_{s+1} D^{s+1} = n_{s+1} D^{s+1}. \]

Since $\CY$ is an approximation cycle, by Definition $h(\CY) \leq a_s D^s$,
and Proposition \ref{bezexthe}.1 implies
\[ h(\CY . \di f) \leq 2\bar{n}_{r+1} D h(\CY) +
   \deg Y 6a \bar{n}_{r+1} D + d \bar{n}_{r+1} D \deg Y \leq \]
\begin{equation} \label{zwischen2}
(\bar{n}_{r+1} a_s + 6 a \bar{n}_{r+1} n_s+ d n_s \bar{n}_{r+1}) D^s \leq 
a_{s+1} D^{s+1}, 
\end{equation}
and for $a>d$ we get
\begin{equation} \label{zwischen3}
t_a(\CY . \di f) \leq (7a + 7) \bar{n}_{s+1}D \deg Y + 
2\bar{n}_{r+1} D h(\CY) < 8 \bar{n}_{s+1} t_a(\CY)
\end{equation}

Further, every irreducible component $\CZ$ of $\CY . \di f$ is a successive
$D$-intersection, hence by Lemma \ref{dh},
\begin{equation} \label{zwischen1}
t_a(\CZ) m_{s+1} D^{t-s} < t_a(\CY) m_s D^{t+1-s}.
\end{equation}

Finally $f_Y^\bot = \bar{f}_Y^\bot$ implies that the restriction of $f$
to $Y$ equals the restriction of $\bar{f}$ to $Y$, and hence
$(\CY. \di f)_\C = Y . \di \bar{f}$.

We make several case distinctions:

{\sc Case 1: $|\di \bar{f},\theta| \leq |Y,\theta|$.} 
By Proposition \ref{bezexthe}.3,
\[ D(\di f . Y, \theta) =
   D(\di \bar{f} . Y, \theta) \leq D(Y,\theta) + 2\bar{n}_{s+1} D h(\CY) + 
   \deg Y 6 a \bar{n}_{s+1} D + d D \deg Y < \]
\[ - b_s t_a(\CY) D^{t+1-s} + 8 \bar{n}_{s+1} t_a(\CY)D, \]
for $a > d$. Since $s \leq t-1$, we have $t+1-s \geq 2$, hence for
big enough $D$, (\ref{zwischen3}) implies
\[ D(\di f . Y, \theta) = D(\di \bar{f} . Y, \theta) \leq
   - \frac{b_s}2 t_a(\CY) D^{t+1-s} \leq - \frac{b_s}{16 \bar{n}_{r+1}} 
   t_a(\di f . \CY) D^{t-s} \leq \]
\[ -b_{s+1} t_a(\di f . \CY) D^{t-s}, \]
by the choice of $b_{s+1}$, and $\bar{n}_{r+1} < \bar{n}_{s+1}$.
By Lemma \ref{appcyc}, $\di f . \CY$ has an irreducible component $\CZ$ with
$D(Z,\theta) \leq - b_{s+1} t_a(\CZ) D^{t-s}$.
Because of (\ref{zwischen2}), the estimate on the degree of $Z$, and
the inequlities $\deg Z \leq \deg (Y . \di f)$, and 
$h(\CZ) \leq h(\CY . \di f)$, $\CZ$ is thus
an approximation cycle, and because of (\ref{zwischen1}), merging
$\CZ$ at the end of $\CC$ one obtains an approximation chain of the form 
$\CC_2$ with $\CC_2 \prec \CC$.

{\sc Case 2:$|\di f^\bot,\theta| \geq |Y,\theta|$.} 

{\sc Case 2.a: $t_a(\CY) \leq 2 \bar{b}_r t_a(\CX_{min}) D^{s-r}/b_s$.}
This inequality together with Proposition \ref{bezexthe}.3, implies
\[ D(\di f . Y, \theta) \leq - \bar{b}_r t_a(\CX_{min}) D^{t+1-r} + 
   2 \bar{n}_{s+1} D h(\CY) + 
   \deg Y 6 a \bar{n}_{s+1} D + d D \deg Y \leq \]
\[ -\frac{b_s}2 t_a(\CY) D^{t+1-s} + 7 \bar{n}_{s+1} t_a(\CY)D. \]
Repeating the argument of case 1, one obtains again
a chain of type $\CC_2$ with $\CC_2 \prec \CC_1$.

{\sc Case 2.b: $t_a(\CY) > 2\bar{b}_r t_a(\CX_{min}) D^{s-r}/b_s$.}
In this case firstly, since 
$m_r \leq \frac{2 \bar{b}_r m_sa_r}{b_{s+1}a_s}$,
by the choice of $b_s$,
\[ t_a(\CX_{min}) m_r D^{t+1-r} \leq
   t_a(\CX_{min}) \frac{2 \bar{b}_r m_sa_r}{b_{s+1}a-s} D^{t+1-r} <
   t_a(\CY) \frac{m_s a_r}{a_s} D^{t+1-s} \leq a_r m_s D^{t+1}. \]
Thus, 
\[ t_a(\CX_{min}) \leq a_r \frac{m_s}{m_r} D^r \leq a_r D^r, \]
and by Lemma \ref{hoch}, either $\CX_{min}$ is an approximation cycle, in
which case $\CX_{min}$ being a component of a locally complete
$D$-intersection would be the end term of a $(D,a)$-approximation chain
\[ \CC_1: \quad \CY_0 \supset \cdots \supset \CX \]
with $\CC_1 \prec \CC$, or $\di f. Y$ is an approximation cycle, and
thus as in case 1
contains an irreducible approximation cycle $\CZ$ gives rise to
a chain of kind $\CC_2$ with $\CC_2 \prec \CC$.

\proof {\sc of Theorem \ref{appcycle}}
Let $r<t$, and 
\[ \CC: \quad \CY_0 \supset \cdots \supset \CY_r, \]
be a chain in $C_{t,D}$ with corresponding global sections $f_1, \cdots f_s$.
It has to be shown that if $s \neq t$, there is a $\CC_2$ with
$\CC_2 \prec \CC$. 

By the previous Lemma there either is a $(D,a)$-approximation chain
\[ \CC_1: \quad \CY_0 \supset \cdots \supset \CX, \]
with $\CC_1 < \CC$ or a $(D,a)$-approximation chain
\[ \CC_2:\quad \CY_0 \supset \cdots \CY_r \supset \CY_{r+1} \]
with $\CC_2 < \CC$. The Theorem follows.

\satz{Corollary of the proof} \label{cor1}
For every $D \in \N$, and every $s \leq t$ there is an approximation cycle
$\CY_{D,s}$ of order $D$ and codimension $s$.
\end{Satz}

\proof
By the proof, starting with an approximation cycle $\CY_1$ of codimension one,
there is a series of approximation cycles
\begin{equation} \label{belcod}
\CY_1, \CY_2, \ldots, \CY_{k-1}, \CY_k
\end{equation}
such that $\CY_k$ has codimension $t$ and for every $l=1, \ldots, k-1$
we either have that $\CY_{l+1}$ is a subvariety of codimension one in 
$\CY_l$, or $\CY_l$ is a subvariety of $\CY_{l+1}$. Hence, (\ref{belcod})
contains at least one cycle of every codimension between $1$ and $t$.

\subsection{Algebraic points with small distance}

Define the sets
\begin{equation} \label{RD}
R_D := \left\{ \alpha \in \Pe^t(\bar{\Q}) | 
     \begin{array}{l} \mbox{$\alpha$ is
                the last term in a} \\ \mbox{successive $(D,a)$-intersection}
     \end{array} \right\}. 
\end{equation}
for every $D \in \N$.
For a generic point 
$\theta \in \Pe^t(\C)$ choose a $\beta_D \in R_D$ which is the end term
of a $(D,aD)$-approximation chain which is minimal with respect to
the relation. By Theorem \ref{appcycle} such a $\beta_D$ exists for every
$D \in \N$, and
by Lemma \ref{appcyc}, and Theorem \ref{appcycle},
$\log |\beta_D,\theta|\leq \varphi_{a,\theta}(\beta_D,\theta) \leq -b_r D$, 
thus
\[ \lim_{D \to \infty} \beta_D = \theta, \]
and 
\begin{equation} \label{N}
\bar{M}:= 
\{ D \in \N|\varphi_{a,\theta}(\beta_{D+1})< \varphi_{a,\theta}(\beta_D)\} 
\end{equation}
is an infinite set.

We will still use the constants (\ref{const}) and further choose a
positive number $\bar{b}$ that is sufficiently small
compared with the constants (\ref{const}).
(How small it has to be will become clear within the following proofs,
i.\@ e.\@ within the following proofs there will appear finitely many 
positive expressions in terms of 
the constants (\ref{const}) all of which are needed to be
bigger than $\bar{b}$.) Further define 
\begin{equation} \label{const2}
n:= \left[\frac{8\bar{n}_t}{b_t}\right]+1,
\end{equation}

\satz{Definition}
Let $a >>0$, $D \in \N$, and $\bar{D} = nD$.
A triple $(f,\bar{f},\CY)$ consisting of a global sections 
$f \in \Gamma(\Pe^t,O(D))_\Z, \bar{f} \in \Gamma(\Pe^t,O(D))_\C$
and a $\Z$-irreducible subvariety $\CY \subset \Pe^t_\Z$ of codimension $t$
with $t_a(\CY) \leq (a_t+an_t) \bar{D}^t$
is called a $(D,a)$-approximation triple if 
$f_Y^\bot = \bar{f}_Y^\bot \neq 0$, 
hence $f|_Y = \bar{f}|_Y \neq 0$, and
\[ D(Y,\theta) \leq -b_t t_a(\CY) \bar{D}, \quad
   \log |f_\CY^\bot| \leq 6 a \bar{n}_t D, \quad
   \log |\la \bar{f} | \theta \ra| \leq - \frac{\bar{b}}{n^{6t}} t_a(\CY) D, \]
\end{Satz}

\satz{Proposition} \label{teil1}
If $(f,\bar{f},\CY)$ is a $(D,a)$-approximation triple, then
\[ \log |Y,\theta| \leq - b t_a(\CY) \bar{D}, \]
with a constant $b>0$ only depending on $t$.
\end{Satz}

\satz{Proposition} \label{teil2}
There is an infinite subset $M \subset \N$ such that for each $D \in M$
there exists a $(D,a)$-approximation triple $(f,\bar{f},\CY)$ with
\[ t_a(\CY) \leq k a \bar{D}^t, \quad \mbox{hence} \quad
   \sqrt[t+1]{\frac{t_a(\CY)}ka} \leq \bar{D}^{\frac t{t+1}}, \]
with a number $k \in \N$ only depending on $t$.
\end{Satz}

Lets first see why the two Propositions together imply Theorem \ref{main}
for the case of projective space.
For $D \in $M from Proposition \ref{teil2}, define
\[ D_1 := 
   \left[\left( \frac{t_a(\CY) \bar{D}}a 
   \right)^{\frac1{t+1}}\right]. \]
Then, $D_1 \geq \bar{D}^{\frac1{t+1}}$, and because of the estimate
on $t_a(\CY)$,
\[ \deg Y \leq \frac{t_a(\CY)}a \leq D_1^t, \quad
   h(\CY) \leq t_a(\CY_{\bar{D}}) \leq a D_1^t. \]
Further, by Proposition \ref{teil1},
\[ \log |Y,\theta| \leq - b t_a(\CY) \bar{D} 
   = -b a D_1^{t+1}, \]
and Theorem \ref{main} follows.

\subsubsection{Using approximation triples}

I present two different proofs of Propostion \ref{teil1} one using
Proposition \ref{bezexthe}.2, the other Proposition \ref{bezmult}. Of course
these proofs does not deliver two independent proofs of the main Theorem, 
because the two Propositions were obtained essentially by a single proof.

\proof 1:
Let $Z := \di \bar{f}$. 
Firstly $|Z,\theta| \geq |Y,\theta|$, since otherwise
Proposition \ref{bezexthe}.3 for $a \geq d$ would imply
\[ 0 = D(\di f . Y,\theta) \leq
   D(Y,\theta) + 2D h(\CY) + \deg Y  \log |f_Y^\bot| + d D \deg Y \leq \]
\[ - b_t t_a(\CY) n D + 7 \bar{n}_t t_a(\CY) D, \]
which because of $n = \bar{D}/D \geq 8\bar{n}_t/b_t$ would be a contradiction.

Next, by Proposition \ref{bezexthe},2 for $\nu, \kappa$ defined by $f_{Z,Y}$,
\begin{eqnarray*}
\nu \kappa \log |\theta,Y| &=& \nu \kappa \log |\theta,Y + Z| \\ \\
&\leq &\kappa D(\theta,Y) + \nu D \log |\la \bar{f}|\theta\ra| +
   6 a \bar{n}_t D \deg Y + 2D h(\CY) + d_1 D \deg Y \\ \\
&\leq& -\kappa b_t t_a(\CY) \bar{D} - \nu \bar{b}/n^{6t}  t_a(\CY) D +
  7 D \bar{n}_t t_a(\CY). 
\end{eqnarray*}
Assume first that 
$\deg Z + \deg Y > \frac{8\bar{n}_t}{\min(\bar{b}/n^{6t},b_t)}$.
Then, choosing $t \in [0,1]$ in Theorem \ref{bezout}.3 in such a way that
$\nu + \kappa = [\frac{8\bar{n}_t}{\min(\bar{b}/n^{6t},b_t)}]$,  we get
\[ \nu \kappa \log |\theta,Y| \leq (\nu + \kappa)
   \min(\bar{b}/n^{6t},b_t) t_a(\CY) D + 7 D t_a (\CY) \leq \]
\[ -8 \bar{n}_t t_a(\CY) D + 7 D \bar{n}_t t_a(\CY) \leq 
   -\bar{n}_tt_a(\CY) D. \]
As $\nu \kappa \leq \frac{(\nu + \kappa)^2}4 = 
\frac{16\bar{n}_t^2}{\min(\bar{b}/n^{6t},b_t)^2}$, 
this in turn implies
\[ \log |\theta,Y| \leq 
   -\frac{\min(\bar{b}/n^{6t},b_t)^2}{16n^{6t}\bar{n}_t} t_a(\CY) \bar{D}, \]
proving the claim with 
$b = \frac{\min(\bar{b}/n^{6t}/n^{6t},b_t)^2}{16n^{6t}\bar{n}_t}$.
If, on the other hand 
$\deg Z + \deg Y \leq \frac{8 \bar{n}_t}{\min(\bar{b}/n^{6t},b_t)}$, 
then by Theorem 
\ref{bezout}.1,
\[ \log |Y,\theta| \leq \frac1{\deg Y} \; 
   D(Y,\theta) + c \leq 
   -\frac{\min(\bar{b}/n^{6t},b_t)}{8 \bar{n}_t} b_t t_a(\CY) \bar{D}, \]
proving the claim with $b=\frac{b_t min(\bar{b}/n^{6t},b_t)}{8 \bar{n}_t}$.

\vspace{2mm}

\proof 2:
We have again $|Y,\theta|\leq |Z,\theta|$ with $Z = \di f$ as in the first 
proof. Let 
\[ Y_\C = \sum_{i = 1}^{\deg Y} y_i, \]
with 
\[ |y_1,\theta| \leq |y_2,\theta| \leq \cdots \leq |y_{\deg Y}, \theta|, \]
and points counted with multiplicities.
If $\deg Y \leq m_3 := [8n^{6t}\bar{n}_t/\bar{b}]+1$, with $n$ from 
(\ref{const2}), clearly $\log |Y,\theta| \leq \frac1{\deg Y} D(Y,\theta) \leq
- \frac{b_t}{m_3} t_a(\CY) \bar{D}$, proving the claim with
$b = b_t/m_3$.

If on the other hand $\deg Y \geq m_3$,
I claim, that $|y_{m_3},\theta| \geq |X,\theta|$.
Assume the opposite; then, by Proposition \ref{bezexthe}.2,
with $t$ such that $f_{X,Y}(t) = (0,m_3)$,
\[ 0 = D(\theta, Y. Z) \leq
   m_3 \log |\la \bar{f}|\theta\ra| + 
   \deg Y 6 a D + 2 D h(\CY) + d D \deg Y \leq \]
\[ - m_3 \frac{\bar{b}}{n^{6t}} t_a(\CY) D + 7 \bar{n}_t t_a(\CY) D. \]
which because of $m_3 \geq \frac{8n^{6t}\bar{n}_t}{\bar{b}}$ is a 
contradiction.

Further, by Proposition \ref{bezmult},
\begin{equation} \label{bew2eins}
\sum_{i=m_3+1}^{\deg Y} \log |y_i,\theta| \geq - 2 D h(\CY) -
\deg Y \log |f_Y^\bot|_{L^2(\Pe^t)} - d \deg X \deg Y \geq 
- 7 \bar{n}_tD t_a(\CY), 
\end{equation}
and
\begin{equation} \label{bew2zwei}
\sum_{i=1}^{\deg Y} |\theta,y_i|=D_{pt}(Y,\theta) \leq - b_t t_a(\CY) \bar{D}, 
\end{equation}
by assumption.
Now (\ref{bew2eins}), and (\ref{bew2zwei}) together imply
\begin{eqnarray*} 
\sum_{i=1}^{m_3} \log |\theta,y_i| &=&
   \sum_{i=1}^{\deg Y} \log |\theta,y_i|-
   \sum_{i=m_3+1}^{\deg Y}\log|\theta,y_i| \\ \\
&\leq& -b_t t_a(\CY) \bar{D} + 7 D t_a(\CY) \\ 
&=& -b_t t_a(\CY) n D + 7 D \bar{n}_t t_a(\CY),
\end{eqnarray*} 
which because of $n \geq 8\bar{n}_t/b_t$ is less or equal
$-\frac1n t_a(\CY) \bar{D}$. Hence,
\[ \log |Y,\theta| = \log |y_1, \theta| \leq
   \frac1{m_3}\sum_{i=1}^{m_3} \log |\theta,y_i| \leq
   - \frac{1}{nm_3} t_a(\CY) \bar{D}, \]
proving the proposition with $b = \frac{1}{nm_3}$.

\subsubsection{Existence of approximation triples}

To proof Proposition \ref{teil2}, we will need to compare approximation
chains $\CC_1 \in C_{D,t}, \CC_2 \in C_{\bar{D},t}$ for $D \neq \bar{D}$.

Let $D<\bar{D}$, and $\alpha_{\bar{D}}$ be an irreducible
approximation cycle of codimension $t$ forming
the end term of a minimal $(D,a)$-approximation chain of order $\bar{D}$. 
Define
\[ C_{D,t}(\alpha_{\bar{D}}) = \{ \CC: \CY_0 \supset \cdots \supset \CY_l 
   \in C_{D,t}| \alpha_{\bar{D}} \subset \CY_l \}. \]
Clearly, $C_{D,t}(\alpha_{\bar{D}})$ is always nonempty, because it contains
the trivial chain.
On $C_{D,t}(\alpha_{\bar{D}})$ we still have the relation $\prec$, and
again, because of the finiteness of $C_{D,t}(\alpha_{\bar{D}})$ 
the relation $\prec$ restricted to the set $\C_{D,t}(\alpha_{\bar{D}})$ has
at least one minimal element in $C_{D,t}(\alpha_{\bar{D}})$.
These minimal elments now usually do not have length $t$, and if
they are not, because they
are minimal it is possible to construe global sections with small
evaluation at $\theta$ that are nonzero on $\alpha_{\bar{D}}$ if
$\bar{D}$.
The fundamental Lemma for this construction is the following.

\satz{Lemma} \label{teil2l1}
If $\alpha_{\bar{D}}$ is not the end term of any minimal chain 
$\CC \in C_{D,t}$, or otherwise said $C_{D,t}(\alpha_{\bar{D}})$
does not contain a minimal element of length $t$, let
$\CY_r, r \leq t-1$ be the last term of such a minimal chain in
$\CC_1 \in C_{D,t}(\alpha_{\bar{D}})$. Then, there is
a locally complete $D$-intersection $\CX$ at $\alpha_{\bar{D}}$  of
codimension $p \leq r$ that contains $\CY_r$, and global sections 
$g \in \Gamma(\Pe^t,\bar{n}_{p+1} D)_Z, \bar{g} \in \Gamma(\Pe^t,\bar{n}_{p+1} D)_\C$
such that $g_{Y_r}^\bot = \bar{g}_{Y_r}^\bot \neq 0$, and
\[ \log |g_{Y_r}^\bot| \leq 6 a \bar{n}_{p+1} D, \quad 
   \log |\la \bar{g}|\theta\ra| \leq - \bar{b}_p t_a(\CX_{min})D^{\dim \CX}, \]
and the restrictions of $g$, and $\bar{g}$ to every irreducible
component of $X$ are non zero.
Furthermore, $g$ has nonzero restriction to $\alpha_{\bar{D}}$, and
\[ t_a(\CX_{min}) m_p D^{\dim \CX} \geq t_a(\CY_r) m_r D^{\dim \CY_r}. \]
\end{Satz}

\proof
If 
\[ \CC_1: \quad \CY_0 \supset \cdots \supset \CY_r \]
is a minimal chain in $C_{D,t}(\alpha_{\bar{D}})$, Corollary
\ref{cor} implies the existence of $\CX$, $g$, and $\bar{g}$, with
the properties stated in the Lemma. If $g$ had zero 
restriction to $\alpha_{\bar{D}}$, then $\CC_1$ would not be minimal in 
$C_{D,t}(\alpha_{\bar{D}})$, because by the proof of Lemma
\ref{zwischen}, case 1, and case 2.a, $\di g . \CY_r$ would again
be a $(D,a)$-approximation cycle containing $\alpha_{\bar{D}}$ and by Lemma
\ref{appcyc} likewise some irreducble component $\CZ$ of $\di g . \CY_r$ would
be a $(D,a)$-approximation cycle giving rise to a chain
\[ \CC_2: \quad \CY_0 \supset \cdots \supset \CY_r \supset \CZ \]
with $\CC_2 \in C_{D,t}(\alpha_{\bar{D}})$, and $\CC_2 \prec \CC_1$ 
contradicting the minimality of $\CC_1$ in $C_{D,t}(\alpha_{\bar{D}})$.

To prove the last claim on  the dimensional $a$-sizes, assume the opposite. 
Then, by Lemma \ref{zwischen}, either $\CX_{min}$
would be a $(D,a)$-approximation cycle, henc $\CC_3 \prec \CC_1$, 
with $\CC_3$ a chain with end term $\CX_{min}$, and
$\CC_1$ were not minimal in $C_{D,t}(\alpha_{\bar{D}})$,  or again, some
$\CZ \subset \di g . \CY_r$ would be a $(D,a)$-approximation cycle containing 
$\alpha_{\bar{D}}$ in which
case again $\CC_1$ were not minimal in $C_{D,t}(\alpha_{\bar{D}})$.

\proof {\sc of Proposition \ref{teil2}}
Recall the definition of the sets $R_D$ and $\bar{M}$ in (\ref{RD}).
We will prove that for each $D \in \bar{M}$, there is a $\bar{D} >D$, and
a $(\bar{D},a)$ approximation triple $(f,\bar{f},\CY)$.
By Theorem \ref{appcycle}, for each $D \in \N$, one can
choose a $(D,a)$ approximation cycle $\alpha_D$ which is the end term of an
minimal $(D,a)$-approximation chain, and has minimal weighted algebraic
distance among all end terms of minimal $(D,a)$-approximation
chains. In particular
\[ \deg \alpha_D \leq n_t D^t, \quad h(\alpha_D) \leq a_t D^t, \quad
   \varphi_{a,\theta} (\alpha_D) \leq - b_t D. \]
With $D \in \bar{M}$ we make the case distinction as to whether
$\varphi_{a,\theta}(\alpha_{\bar{D}}) > -n^2 b_t \bar{D}$ for every
$\bar{D}>D$ or
$\varphi_{a,\theta}(\alpha_{\bar{D}}) \leq -n^2 b_t \bar{D}$ for some 
$\bar{D} > D$, with $n$ the constant from (\ref{const2}).

\vspace{2mm}

{\sc Case 1: $\varphi_{a,\theta}(\alpha_{\bar{D}}) \leq -c_1 b_t \bar{D}$}
for some $\bar{D} > D$.

\satz{Lemma}
In Case 1, for $\bar{D}$ the smallest number $>D$ such that
$\varphi_{a,\theta}(\alpha_{\bar{D}}) \leq -c_1 b_t \bar{D}$, 
the variety $\alpha_{\bar{D}}$ 
is not the end term of any $(\bar{D}-1,a)$-approximation chain, i.\@ e.\@
the minimal elements in $C_{\bar{D}-1,t}(\alpha_{\bar{D}})$ have length 
$\leq t-1$.
\end{Satz}

\proof
A: $\bar{D} > D + 1$. In this case, by the choice of $\bar{D}$,
\[ \varphi_{a,\theta} (\alpha_{\bar{D}-1}) > - c_1 b_t (\bar{D}-1) >
   - c_1 b_t \bar{D} \geq \varphi_{a,\theta}(\alpha_{\bar{D}}), \]
in particular $\alpha_{\bar{D}} \neq \alpha_{\bar{D}-1}$ by Lemma 
\ref{appcyc}, and
hence $\alpha_{\bar{D}}$ is not the end term of a minimal 
$(\bar{D}-1,a)$-approximation chain by the choce of $\alpha_{\bar{D}-k}$
as a minimal element in $C_{D,t}$ that has minimal weighted algebraic
distance to $\theta$. 

B: $\bar{D} = D+1$. In this case, 
$\alpha_{D+1}$ is not the end term of a minimal $(D,a)$-approximation 
chain, since $D \in \bar{M}$.

\vspace{2mm}

From on now let $D = \bar{D}-1$ with $\bar{D}$ from the lemma. This number is not 
necessarily contained in $\bar{M}$ but making this replacement for every $D$
defines another infinite set $\bar{M}_1 \subset \N$.

For $D,\bar{D},\alpha_{\bar{D}}$ as from the previous Lemma let
$\CY_r, \CX, g, \bar{g}$ be an irreducble variety, a locally complet
$D$-intersecting and global section with the properties of Lemma \ref{teil2l1}.
Define the subset 
$C_{\bar{D},t}(\alpha_{\bar{D}},\CY_r)$ consisting
of the succesive $\bar{D}$-intersection chains with end term
$\CY_s$ an irreducible subvariety of codimesnion $s$ fullfilling
\[ \alpha_{\bar{D}} \subset \CY_s \subset \CY_r, \quad
   h(\CY_s) \leq a_s \bar{D}_s, \quad D(\CY_s,\theta) 
   \leq -b_s t_a(\CY) D^{t+1-s}. \]
Note that $\CY_s$ need neither be a $(D,a)$- nor a $(\bar{D},a)$-approximation
cycle. Its purpose is to give a good estimate on $t_a(\alpha_{\bar{D}})$.
Nonetheless on $C_{\bar{D},t}(\alpha_{\bar{D}},\CY_r)$ one has the
relation $\prec$ being defined in the same way as on $C_{\bar{D},t}$ with the
condition on dimensional heigth of order $\bar{D}$ as second condition.

Observe, that although
$\CY_r$ belongs to a minimal chain in $C_{D,t}(\alpha_{\bar{D}})$ with
respect to $\prec$, it need not belong to a minimal chain in
$C_{\bar{D},t}(\alpha_{\bar{D}},\CY_r)$ because $\bar{D}>D$.
Nonetheless, $C_{\bar{D},t}(\alpha_{\bar{D}},\CY_r)$
is again finite, and therby has minimal elments.

\satz{Lemma} \label{teil2l2}
With the above notation, for every $\CY_s$ the end term of a minimal chain 
in $C_{\bar{D},t}(\alpha_{\bar{D}},\CY_r)$, we have 
$s:= codim \CY_s \geq r$. Further there exists a $\CY_s$ which is
the end term of a minimal chain 
$\CC_2 \in C_{\bar{D},t}(\alpha_{\bar{D}},\CY_r)$ that fullfills
\[ m_s t_a(\CY_s) \leq t_a(\CY_r) m_r \bar{D}^{s-r}. \]
For this $\CY_s$ 
there is a locally complete $(\bar{D},a)$-intersection $\bar{\CX}$
of codimension $q \leq s$,
and $f \in \Gamma(\Pe^t,O(\bar{n}_{s+1} \bar{D}))_\Z$, 
$\bar{f} \in \Gamma(\Pe^t,O(\bar{n}_{s+1} \bar{D}))_\C$ with
$f_{Y_s}^\bot = \bar{f}_{Y_s}^\bot \neq 0$, and
\[ \log |f_{\CY_s}^\bot| \leq 6 a \bar{n}_{q+1} \bar{D}, \quad
   \log |\la \bar{f}|\theta\ra| \leq 
   - b_s t_a(\bar{\CX}_{min}) \bar{D}^{\dim \CX}, \]
such that the restrictions of $f,\bar{f}$ to 
every irreducible component of $\bar{\CX}$ are nonzero. If 
furter $s \leq t-1$, then the restriciont of $f,\bar{f}$ to
$\alpha_{\bar{D}}$ are also non zero.
\end{Satz}

\proof
The inequality $s \geq r$ for any $\CY_s$ as in the Lemma is trivial.

Further, since cleary $\CY_r$ is the end term of a chain 
$\CC_1 \in C_{\bar{D},t}(\alpha_{\bar{D}},\CY_r)$,
there is a minimal $\CC_2 \in C_{\bar{D},t}(\alpha_{\bar{D}},\CY_r)$, with
$\CC_2 \prec \CC_1$, hence the end term $\CY_s$ of $\CC_s$ fullfills
\[ t_a(\CY_s) m_s \bar{D}^{t+1-s} < t_a(\CY_r) m_r \bar{D}^{t+1-r}, \]
and the on the dimensional size of $\CY_s$ follows.

By corollary \ref{cor}, there is a locally complete $(\bar{D},a)$-intersecion
$\bar{\CX}$ at $\CY_s$, and global sections 
$f \in \Gamma(\Pe^t,O(\bar{n}_{s+1} \bar{D}))_\Z$,
$\bar{f} \in \Gamma(\Pe^t,O(\bar{n}_{s+1} \bar{D}))_\C$
fullfilling the equalities and inequalities 
in the Lemma, and having nonzero restriction to $\CY_s$, and
every irreducible component of $\bar{\CX}$.
It remains to be proved that $f$ has nonzero restriction to $\alpha_{\bar{D}}$ 
if $s\leq t-1$. Assume $f$ is zero on $\alpha_{\bar{D}}$. Then
$\di f \cap \CY_s$ contains $\alpha_{\bar{D}}$, and from the Theorem of
B\'ezout, and Proposition
\ref{bezexthe}.1, one could deduce a bound on the height and degree 
$\di f \cap \CY_s$, and an using a suitable irreducible cycle $\CZ$ in
$\di f \cap \CY_s$ containing $\alpha_{\bar{D}}$, Lemma \ref{zwischen}
would again supply a chain $\CC_3$ with $\CC_3 \prec \CC_2$ contradicitng
the minimality of $\CC_2$.

\satz{Lemma} \label{teil2l3}
With $D,\bar{D},\alpha_{\bar{D}}$ as above, let $\CY_r,\CX$, $g,\bar{g}$
be the subvarieties and global sections existing by Lemma \ref{teil2l1}, and
$\CY_s,\bar{\CX},f, \bar{f}$ 
the varieties and global section existing by the previous Lemma. Then,
\[ t_a(\alpha_{\bar{D}}) m_t \bar{D} \leq
   2 \frac{m_t}{m_s}\mbox{max} \left(\frac{\bar{b}_s}{b_t}m_q,\bar{b}_q\right) 
   t_a(\bar{\CX}_{min}) D^{\dim \bar{\CX}}, \]
and
\[ t_a(\alpha_{\bar{D}}) m_t \bar{D} \leq
   2\frac{m_t}{b_t} \mbox{max}\left(\bar{b}_s,\frac{m_s}{m_q} \bar{b}_q \right)
   t_a(\CY_s) \bar{D}^{\dim \CY_s}. \]
If $\CY_s = \alpha_{\bar{D}}$, additionally 
\[ \frac12 t_a(\alpha_{\bar{D}}) \bar{D}^{t+1-r} \leq
t_a(\alpha_{\bar{D}}) m_p \bar{D}^{r-s} D^{t+1-r} \leq
   2 m_r \frac{\bar{b}_p}{b_r} \frac{m_p}{m_t} t_a(\CX_{min}) D^{\dim \CX}. \]
\end{Satz}

\proof {\sc of Propostion \ref{teil2}, Case 1, finish:}
With the notations of the two previous Lemmas, assume first
that $\alpha_{\bar{D}}$ is properly contained in $\CY_s$. Then, the 
restrictions of $f,\bar{f}$ to $\alpha_{\bar{D}}$ are non zero, and we get 
$D(\alpha_{\bar{D}},\theta) \leq - b_t n^2 t_a(\alpha_{\bar{D}} \bar{D}$. 
Secondly, by Lemma \ref{teil2l2},
$\log |\bar{f}_{\alpha_{\bar{D}}}^\bot|\leq 
\log |\bar{f}_{Y_s}^\bot|\leq 6 a \bar{n}_{q+1} \bar{D}$.
Finally, again by Lemma \ref{teil2l2},
\[ \log |\la \bar{f}|\theta \ra| \leq - \bar{b}_q t_a(\bar{\CX}_{min}) 
   D^{\dim \bar{\CX}}, \]
which by Lemma \ref{teil2l3} is at most
\[ -\frac12 b_t \frac{m_s}{m_t} 
   \mbox{max} \left(\frac{1}{\bar{b}_q},\frac{\bar{b}}{b_t}m_q\right) 
   t_a(\alpha_{\bar{D}}) \bar{D}, \]
Thus, with 
$\bar{b} \leq \frac12 b_t \frac{m_s}{m_t}
\left(\frac{1}{\bar{b}_q},\frac{\bar{b}}{b_t}m_q\right)$, we
have the three inequalities
\[ D(\alpha_{\bar{D}},\theta) \leq - b_t n^2 t_a(\alpha_{\bar{D}}) \bar{D}, 
   \quad
   \log |f_{\alpha_{\bar{D}}}^\bot| \leq 6a \bar{n}_{q+1} \bar{D} \leq
   6a \bar{n}_t n \bar{D}, \]
\[ \log |\la\bar{f}|\theta\ra| \leq - \bar{b} t_a(\alpha_{\bar{D}}) \bar{D} 
   \leq - \frac{\bar{b}}{n^{6t}} t_a(\alpha_{\bar{D}}) n \bar{D}, \]
which are the 3 inequalities in the Definition of approximation triples
with $\bar{D}$ replaced by $n \bar{D}$. Since
\[ t_a(\alpha_{\bar{D}}) \leq a \deg \alpha_{\bar{D}} + h(\alpha_{\bar{D}}) 
   \leq (a n_t+a_t) \bar{D}^t \leq (a n_t + a_t) (n\bar{D})^t, \]
this proves that $(f, \bar{f},\alpha_{\bar{D}})$ is an approximation triple
of ordet $n \bar{D}$.

\vspace{2mm}

If on the other hand $\alpha_{\bar{D}} = \CY_s$, then firstly still
$D(\alpha_{\bar{D}}, \theta) \leq - b_t n^2 t_a(\alpha_{\bar{D}})\bar{D}$, 
secondly, by Lemma \ref{teil2l1}, the restriction
of $g,\bar{g}$ to $\alpha_{\bar{D}}$ is non zero, and
$\log |g_{\alpha_{\bar{D}}}^\bot| \leq |g_{Y_r}^\bot| \leq 6 a \bar{n}_{p+1} D$,
and thirdly also by Lemma \ref{teil2l1},
\[ \log |\la \bar{g}|\theta\ra| \leq - \bar{b}_p t_a(\CX_{min}) D^{t+1-p}, \]
which by Lemma \ref{teil2l3}, and the face $\bar{D}=D+1$ is at most
\[-\frac14 \frac{m_t}{m_r} \bar{b}_r  t_a(\alpha_{\bar{D}}) \bar{D}, \]
and with $\bar{b} \leq -\frac14 \frac{m_t}{m_r} \bar{b}_r$ we have again
three inequalities
\[ D(\alpha_{\bar{D}},\theta) \leq -b_t n^2 t_a(\alpha_{\bar{D}})\bar{D}, \quad
   \log |g_{\alpha_{\bar{D}}}^\bot| \leq 6 a \bar{n}_t n \bar{D}, \quad
   \log |\la \bar{g} |\theta\ra| \leq - \frac{\bar{b}}{n^{6t}} 
   t_a(\alpha_{\bar{D}}) n \bar{D}, \]
and $(g,\bar{g},\alpha_{\bar{D}})$ is an approximation triple of
order $n \bar{D}$,
since again $t_a(\alpha_{\bar{D}}) \leq (a_t+an_t) \bar{D}^t$.

\proof {\sc of Lemma \ref{teil2l3}:}
Let $\CY_s,\bar{\CX}$, $f,\bar{f}$ be the subvarieties and global sections
from Lemma \ref{teil2l3}, and
\[ \CC: \quad \CY_0 \supset \cdots \supset \alpha_{\bar{D}} \]
the minimal element of $C_{\bar{D},t}$ to which $\alpha_{\bar{D}}$
belongs.
If $|\di \bar{f},\theta|$ were smaller than $|\alpha_{\bar{D}},\theta|$, 
Proposition \ref{bezexthe}.3 would imply
\[ 0 = D(\di \bar{f} . \alpha_{\bar{D}},\theta) \leq
   D(\alpha_{\bar{D}},\theta) + 7 \bar{n}_{q+1} \bar{D}t_a(\CY) \leq 
   -b_t n^2 t_a(\CY) \bar{D} + 7 \bar{n}_{q+1} \bar{D} t_a(\CY), \]
which because of $n \geq (8 \bar{n}_t/b_t\leq(8 \bar{n}_{q+1})/b_t$ 
is impossible. Hence, 
$|\di \bar{f},\theta| \geq |\alpha_{\bar{D}},\theta|$.
Further $\di f$ intersects $\CY_s, \CY_r$, and each component of 
$\bar{\CX},\CX$ properly. The first by Lemma \ref{teil2l2}, the second because
$\CY_r$ contains $\CY_s$, and is irreducible, and the third by construction.
Finally, the proper intersection with every component of $\CX$ follows from 
the fact that $\CX$ is a locally complete intersection at $\CY_r$, and
$\di f$ already intersects $\CY_s$ and thereby $\CY_r$ properly.

Now, if $2t_a(\bar{\CX}_{min}) \bar{b}_q \bar{D}^{\dim \bar{\CX}} <
   t_a(\alpha_{\bar{D}}) b_t \bar{D}$, since
$D(\di \bar{f}. \alpha_{\bar{D}},\theta) = 0$, the proof of
Lemma \ref{hoch} implies that
$\bar{\CX}_{min}$ is a $(\bar{D},a)$-approximation cycle. Being a locally
complete $\bar{D}$-intersection, it belongs to a 
$(\bar{D},a)$-approximation chain
$\CC_1$. Since $\bar{\CX}_{min}$ contains $\alpha_{\bar{D}}$, we get
$\CC_1 \prec \CC$ contradicting the minimality of $\CC_1$ in 
$C_{\bar{D},t}$ this is a contradiction. Thus,
\[ 2t_a(\bar{\CX}_{min}) \bar{b}_q \bar{D}^{\dim \bar{\CX}} \geq
   t_a(\alpha_{\bar{D}}) b_t \bar{D}, \]
and hence
\begin{equation} \label{schluss1}
2t_a(\bar{\CX}_{min}) \frac{m_t}{b_t}\bar{b}_q \bar{D}^{\dim \bar{\CX}} \geq
   t_a(\alpha_{\bar{D}}) m_t \bar{D}. 
\end{equation}

{\sc Case 1: $t_a(\bar{\CX}_{min}) m_q \bar{D}^{\dim \bar{\CX}} \leq
t_a(\CY_s) m_s \bar{D}^{\dim \CY_s}$.} In this case from (\ref{schluss1})
the inequality
\begin{equation} \label{l3eins}
t_a(\alpha_{\bar{D}}) m_t \bar{D} \leq 
2t_a(\bar{\CX}_{min}) \frac{m_t}{b_t}\bar{b}_q \bar{D}^{\dim \bar{\CX}} \leq
2 t_a(\CY_s) 
\frac{m_s}{m_q} \frac{m_t}{b_t}\bar{b}_q\bar{D}^{\dim \CY_s} 
\end{equation}
follows.

{\sc Case 2: $t_a(\bar{\CX}_{min}) m_q \bar{D}^{\dim \bar{\CX}} \geq
t_a(\CY_s) m_s \bar{D}^{\dim \CY_s}$.} 

If also 
$t_a(\alpha_{\bar{D}}) b_t \bar{D} \geq 2t_a(\CY_s) \bar{b}_s 
\bar{D}^{\dim \CY_s}$,
again because of $D(\di \bar{f} . \alpha_{\bar{D}},\theta)=0$, the proof of
Lemma \ref{hoch} implies that
$\CY_s$ is an $(\bar{D},a)$-approximation cycle. Beeing a successive
$\bar{D}$-intersection, it belongs to a $(\bar{D},a)$-approximation
chain $\CC_2$ and $\CY_s$ contains $\alpha_{\bar{D}}$, which is 
impossible because $\CC$ is minimal in $C_{\bar{D},t}$. Hence,
\[ t_a(\alpha_{\bar{D}}) b_t \bar{D} \leq 2t_a(\CY_s) \bar{b}_s 
   \bar{D}^{\dim \CY_s}, \]
from which together with the assumption
\begin{equation} \label{l3zwei}
t_a(\alpha_{\bar{D}}) m_t \bar{D} \leq t_a(\CY_s) \frac{m_t}{b_t} \bar{b}_s 
\bar{D}^{\dim \CY_s} \leq
2 \frac{m_t}{b_t} \frac{\bar{b}_s}{m_s} m_q t_a(\CX_{min}) 
\bar{D}^{\dim \bar{\CX}} 
\end{equation}
follows.


Finally, if $\CY_s = \alpha_{\bar{D}}$ assume
$2t_a(\CX_{min}) \bar{b}_p D^{\dim \CX} < t_a(\CY_r) b_r D^{\dim \CY_r}$.
Then Lemma \ref{hoch} implies that either
$\di g. \CY_r$ is a $(D,a)$-approximation cycle, being the end term of an
approximation chain  
$\bar{\CC_2} \in C_{D,t}(\alpha_{\bar{D}})$ with 
$\bar{\CC}_2 \prec \bar{\CC}_1$, where $\bar{\CC}_1$ is the chain in 
$C_{D,t}(\alpha_{\bar{D}})$ with end term 
$\CY_r$, or $\CX_{min}$ is the end term of a chain 
$\bar{\CC}_3 \in C_{D,t}(\alpha_{\bar{D}})$ with 
$\bar{\CC}_3 \prec \bar{\CC}_1$.
Both possibilities contradict the minimality of $\bar{\CC}_1$ 
in $C_{D,t}(\alpha_{\bar{D}})$.
Hence,
\[ 2t_a(\CX_{min}) m_p D^{\dim \CX} \geq t_a(\CY_r) \frac{m_p}{\bar{b}_p} b_r 
   D^{\dim \CY_r} \geq 
  t_a(\CY_s) \frac{m_s}{m_r} \frac{m_p}{\bar{b}_p} b_r 
  \bar{D}^{r-s} D^{t+1-r} ,\]
the last inequality by Lemma \ref{teil2l2}.
Since $\CY_s = \alpha_{\bar{D}}$, hence $s=t$, the claim follows,
if one remembers $\bar{D}=D+1$.

\proof {\sc of Proposition \ref{teil2}, 
Case 2: $\varphi_{a,\theta}(\alpha_{\bar{D}}) \geq -n^2 b_t \bar{D}$
for every $\bar{D} > D$.}

In this case, for every $\bar{D} \geq D$ choose $\beta_{\bar{D}}$
as the end term of a minimal $(\bar{D},a)$-approximation chain which
has minimal dimensional $a$-size
$t_a(\beta_{\bar{D}}) m_t \bar{D}$ among all end terms of minimal
$(\bar{D},a)$- approximation chains. Since $\alpha_{\bar{D}}$ was
chosen with minimal weighted algebraic distance among the end 
terms of minimal elements in $C_{\bar{D},t}$,
\[ \varphi_{a,\theta} (\beta_{\bar{D}})\geq\varphi_{a,\theta}(\alpha_{\bar{D}})
   > - n^2 b_t \bar{D}. \]
Further, by Theorem \ref{appcycle},
$\varphi_{a,\theta} (\beta_{\bar{D}}) \leq -b_t \bar{D}$. 
Hence,
$\varphi_{a,\theta} (\beta_{D_2}) < \varphi_{a,\theta}(\beta_{D,1})$
for $D_1, D_2 > D$, and $D_2 \geq n^2 D_1$, and 
by Lemma \ref{appcyc}, $\beta_{D_2} \neq \beta_{D_1}$.
Assume that
\[ t_a(\beta_{\bar{D}}) \bar{D} < 
   \frac1{n^{6t}} t_a(\beta_{n^2\bar{D}})(n^2\bar{D}) \]
for all $\bar{D} > D$. Then, for every $l \in \N$,
\[ (a_t+ an_t) (n^{2l}D)^{t+1} \geq 
   t_a(\beta_{n^{2l}D}) n^{2l}D > n^{6tl} t_a(\beta_D) D \geq n^{6tl} D, \]
which for big enough $l$ contradicts the fact that$t+1 < 3t$.
Hence, there is some $\bar{D} > D$ such that
\begin{equation} \label{schluss2}
t_a(\beta_{\bar{D}}) \bar{D} \geq \frac1{n^{6t}} 
t_a(\beta_{n^2\bar{D}}) (n^2\bar{D}).
\end{equation}

Let
\[ \CC_1: \quad \CY_0 \supset \cdots \supset \CY_s \]
be a minimal element in $C_{\bar{D},t}(\beta_{n^2\bar{D}})$. 
Then, $s<t$, because otherwise $\beta_{\bar{D}+k}$ would be the end term
of a minimal element in $C_{\bar{D},t}$ contradicting the assumption
$\varphi_{a,\theta} (\beta_{\bar{D}}) > - n^2 b_t \bar{D}$ together with 
the inequality $\varphi_{a,\theta}(\beta_{n^2\bar{D}}) \leq -b_t n^2 \bar{D}$
holding by Theorem \ref{appcycle}.

Let $\bar{\gamma}_{\bar{D}}$ be the end term of a minimal chain
$\CC_3 \in C_{\bar{D},t}$ with $\CC_3 \prec \CC_1$. Then, by the
choice of $\beta_{\bar{D}}$, and the definition of the relation $\prec$,
\[ t_a(\beta_{\bar{D}}) m_t \bar{D} \leq t_a(\gamma_{\bar{D}}) m_t \bar{D} <
   t_a(\CY_s) m_s \bar{D}^{\dim \CY_s}, \]
and by (\ref{schluss2}), 
\begin{equation} \label{teil22}
\frac1{n^{6t}} t_a(\beta_{n^2\bar{D}}) m_t (n^2\bar{D}) \leq
t_a(\beta_{\bar{D}})m_t \bar{D} < t_a(\CY_s) m_s \bar{D}^{\dim \CY_s}.
\end{equation}

By corollary \ref{cor}, there is a
locally complete $\bar{D}$-intersection $\CX$ of codimension 
$r \leq s$ containing $\CY_s$, and $g \in \Gamma(\Pe^t,O(\bar{D}))_\Z$, 
$\bar{g} \in \Gamma(\Pe^t,O(\bar{D}))_\C$ with 
$g_{Y_s}^\bot = \bar{g}_{Y_s}^\bot \neq 0$, and
\[ \log |g_{Y_s}^\bot| \leq 6a \bar{n}_{r+1} \bar{D}, \quad
   \log |\la \bar{g} |\theta\ra| \leq - \bar{b}_s t_a(\CX_{min}) 
   \bar{D}^{\dim \CX}. \]
Because of the minimality of $\CC_1$ in $C_{\bar{D},t}(\beta_{n^2 \bar{D}})$,
and (\ref{teil22}),
\[ -\bar{b}_s t_a (\CX_{min}) \bar{D}^{\dim \CX} \leq
   - \bar{b}_s \frac{m_s}{m_r} t_a(\CY_s) n^2\bar{D}\leq  
   -\bar{b}_s \frac{m_t}{m_r} \frac1{n^{6t}} 
   t_a(\beta_{n^2\bar{D}}) n^2 \bar{D}. \]
Hence, with $\bar{b} \leq \bar{b}_s \frac{m_t}{m_r}$ we get
\[ \log |\la\bar{g}|\theta \ra| \leq - \frac{\bar{b}}{n^{6t}}
   t_a(\beta_{n^2\bar{D}}) n^2 \bar{D}. \]
Together with
\[ \log |g^\bot_{\CY_s}| \leq 6 a \bar{n}_{r+1} \bar{D} \leq
   6a \bar{n}_t \bar{D} \quad \mbox{and} \quad
   D(\beta_{n^2\bar{D}},\theta) \leq -b_t n^2 \bar{D}, \]
and $t_a(\beta_{n^2}\bar{D}) \leq (a_t + a n_t) (n^2 \bar{D})^t$,
this proves that
$(g,\bar{g},\alpha_{\bar{D}+k})$ is an approximation triple with $n$
replaced by $n^2$, and $D$ by $\bar{D}$.

\subsection{Varieties of higher dimension}

To prove that there are also algebriac subvarieties of higher dimension
in $\Pe^t$, that are very close to $\theta$, the general strategy is to
use Corollary \ref{coreinl}, i.\@ e.\@ find subvarieties of bounded height
and degree that contain the subvariety $\alpha_D$ of codimenion $t$
that fullfills Theorem \ref{main}.
However, it is not as easy to find these subvarieties as it was to
find approximation cycles of higher dimension in Corollary \ref{cor1}.
The reason is that although it is possible for the reasons state in
Corollary \ref{cor1} to find approximation cycles of every codimension,
they need not all contain $\alpha_D$ even though $\alpha_D$ belongs
to a minimal element in $\CC \in C_{D,t}$. This is because if in
construeing $\CC \prec \CC_0$ with $\CC_0$ the empty chain only because
in the series $\CC = \CC_1, \CC_2, \ldots, \CC_0 = \CC_l$ from the
definition of $\prec$ there appears an $i$ with
$\CC_i \prec \CC_{i+1}$ such that the end term $\CC_i$ contains the 
end term of $\CC_{i+1}$ one has no guarantee that the end terms 
of $\CC_j, j= 1, \ldots,i$ will contain $\alpha_D$.

However, it can be proved that for generic points $\theta \in \Pe^t$
such that case 2 in the proof of Proposition \ref{teil2} applies, one 
has a series $\CC = \CC_1, \ldots, \CC_t = \CC_0$, and the end term
of $\CC_i$ is a subvariety of codimension $1$ in $\CC_{i+1}$ for
every $k=1,\ldots,t-1$. It can further be proved that these generic points
are almost all points in $\Pe^t$ in the sence of Lebesque measure
(see \cite{Mahler}).
For a generic points $\theta$ for which case 1 in the proof of Proposition 
\ref{teil2} applies, the situation is much more complicated. A proof
that one still finds subvarieties of higher dimension that have small
distance to $\theta$ will also be given in \cite{Mahler}.

\section{The general case}

\proof {\sc of Theorem \ref{main}.2 (arbitrary quasi projective 
scheme over \\ $\spec \; \Z$)}

Firtly, replacing $\CX$ by the algebraic closure of $\{ \theta \}$ it may be
assumed that $\CX=\CY$ is irreducible of relative dimension $t$, and
the algebraic closure of $\{ \theta \}$.
Assume first that $\CL$ is very ample, and choose global sections
$s_1, \cdots, s_{m+1} \in \Gamma(\CX,\CL)$ with $m$ minimal that define an
embedding $i: \CX \to \Pe^n = \Pe(\Z^{n+1})$.
Then $O(1)|_{\CX} = \CL$, 
\[ \deg_{O(1)} i(X) = \deg_L X, \quad
   h_{\overline{O(1)}} (i(\CX)) = h_{\bar{\CL}} (\CX), \]
and for $\alpha \in \CX(\bar{\Q})$,
\[ \deg_{\CL} i(\alpha) = \deg_{\CL}(\alpha), \quad
   h_{\overline{O(1)}} (i(\alpha)) = h{\bar{\CL}}(\alpha). \]
Since we are only interested in algebraic points very close to $\theta$
we may replace $\CX$ by the closure of $i(\CX)$ in $\Pe^m$.
Let $E \cong \Z^{t+1} \subset M,N \cong \Z^{n-t}$ submodules defined by 
any choice of $t+1$ of the $n+1$ coordinates such that $M = E \oplus N$. 
Then the canonical projection $M=E \oplus N \to E$ induces maps
\[ \Pe^n \setminus \Pe(N) \to \Pe(E), 
   \quad \varphi: \CX \setminus \Pe(N) \to \Pe(E). \]
Since $\theta$ is contained in no proper subscheme of $\CX$ it is contained in
$\CX \setminus \Pe(N)$, and $\varphi_\C$ is injective in a certain 
neighbourhood 
$U_\C$ of $\theta$ in the $\C$-topoloty,
because the derivatives of $\varphi$ with respect to the
coordinates of $N$ can't vanish at $\theta$ as they  define algebraic
subvarieties.

Let $V_\C := \varphi(U_\C)$. Then, for any 
$\alpha \in \Pe(E)(\bar{\Q}) \cap U_\C$, the set $\varphi^{-1}(\alpha)$, 
equals the intersection of $\CX(\bar{\Q})$ with the projective subspace 
$\Pe_\alpha \subset \Pe(M)$ corresponding to the submodule
$M_\alpha := N_{\bar{\Q}} \oplus \bar{\Q} \hat{\alpha}$,
where $\hat{\alpha}$ is any vector in $\Pe(E)_{\bar{\Q}}$ representing
$\alpha$. It is easily seen that 
$\deg_{O(1)} \Pe_\alpha = \deg_{O(1)} \alpha$, and
$h_{\overline{O(1)}}(\Pe_\alpha)) = 
h_{\overline{O(1)}}(\alpha) + \sigma_{n-t}$.

Now, $\varphi(\theta)$ is contained in no proper algebraic subset of
$\Pe(E)$ defined over $\bar{\Q}$, and thus
by Theorem \ref{main} for the case of projective space, there is a positive
number $b$ such that for every sufficently big $a$, there is an infinite subset
$\bar{L} \subset \N$, such that for each $D \in \bar{L}$, there is an
$\alpha_D \in V_\C cap \Pe(E)(\bar{\Q})$ with
\[ \deg(\alpha_D) \leq D^t, \quad h(\alpha) \leq aD^t, \quad
   \log |\varphi(\theta), \alpha_D| \leq -b a D^{t+1}. \]

Consider the points
\[ \beta_D := \varphi^{-1}(\alpha_D) = (\Pe_\alpha \cap \CX) \setminus \Pe(M) 
   \in U \cap \CX(\bar{\Q}). \]
Clearly, $\deg \beta_D \leq \deg \CX \deg \alpha_D$, thus
\[ \deg \beta_D \leq \deg X \; D^t. \] 
Further, by the arithmetic B\'ezout Theorem,
\begin{eqnarray*}
 h (\beta_D) &\leq& \deg \Pe_{\alpha_D} h(\CX) + 
   \deg X h (\Pe_{\alpha_D}) + d \deg \Pe_{\alpha_D}  \deg X \\ 
&=& \deg \alpha_D h (\CX) + 
   \deg X (h(\alpha_D) + \sigma_{n-t})+ d \deg \alpha  \deg X \\ 
& \leq & h (\CX) D^t + \deg X a D^t + d \deg X D^t. 
\end{eqnarray*}
Finally, since $U_\C$ may be chosen relatively compact,
there are constants $C,C'> 0$ such that for $x,y \in U_\C$, 
\begin{equation} \label{commet}
C|x,y| \leq |\varphi(x),\varphi(y)| \leq C' |x,y|. 
\end{equation}
Since $\varphi$ is bijective on $U_\C$, for every $\alpha_D$ close 
enough to $\varphi(\theta)$, the $\C$-valued points of $\beta_D$ contain
a point $\bar{\beta}_D \in U_\C \cap X(\bar{Q})$, and
\[ \log |\theta,\beta_D| \leq 
   \log |\theta, \bar{\beta}_D| \leq \log |\varphi(\theta), \alpha_D| 
   - \log C \leq - b a D^{t+1} - \log C. \]
Choosing 
\[ \bar{D} := \left[\frac{D}{\sqrt[t]{\frac{h(\CX)}a + 3/2 \deg X}}\right], \]
and remembering that $a$ can be chosen to be $2\geq d$, we get
\[ \deg \beta_D \leq \bar{D}^t, \quad h(\beta_D) \leq a \bar{D}^t, \quad
   \log |\beta_D,\theta| \leq
   \frac{1}{\left(\frac{h(\CX)}a + 3/2 \deg X\right)^{\frac{t+1}t}} \; 
   b a \bar{D}^{t+1}- \log C. \]
For $D>>0$, replacing $3/2 \deg X$ by $2 \deg X$,
the series $\beta_D$ thus fullfills all requirements of Theorem \ref{main}.

If $\CX$ is only ample, use the property that there are global
sections $s_1, \ldots, s_{m+1}$ of $\CL$ on $\CX$ that define a finite
map $\psi: \CX \to \Pe^m$. By the proof just given the Theorem holds
for $\psi(\CX)$, and thereby also for $\CX$.

\proof {\sc of Theorem \ref{main}.2 (Arbitrary scheme over an
arbitrary number ring)}

To prove the Theorem over an arbitrary ring of integers $\CO_k$ in
a number field $k$ there are two possibilities: The first is to extend the
proofs of 
the metric B\'ezout Theorem in \cite{App1}, the estimates for arithemetic
Hilbert functions in \cite{App2}, and the Propositions in
this paper to projective spaces over $\spec \; \CO_K$ which only entails some 
further effort of a technical kind, and derive the Theorem in the same way.
The second possibility is to use Weyl restriction from
$\spec \; \CO_k$ to $\spec \; \Z$ and use that the Theorem has already been
proved over $\spec \; \Z$. I will here choose the second possibility:

So let $\CO_k$ be the ring of integers of a number field $k$ of degree
$d$ over $k$, let 
$\CX$ be a projective variety of relative dimension $t$ over $\spec \CO_k$,
and $\theta \in \CX(\C_\sigma)$ a generic point, and define
\[ \CX_1 := R_{\CO_k/\Z} \CX. \]
It has relative dimension $t$ over $\spec \; \Z$.

Then, $\CX_1(\C) = \CX(k \otimes_\Q \C)$, 
hence $\theta$ defines a generic point
$\theta_1 \in \CX_1(\C)$. By the Theorem for schemes over
$\spec \; \Z$ there is a $b > 0$ such that for given $a >0$ there is
an infinite subset $M \subset \N$ such that for all $D \in M$ there is
an $\alpha_D \in \CX_1(\bar{\Q})$ with
\[ \deg \alpha_D \leq D^t, \quad h(\alpha_D) \leq a D^t, \quad
   \log |\alpha_D,\theta_1| \leq - \frac1 
   {\left(\frac{h(\CX)}a + 2 \deg X \right)^{\frac{t+1}t}} b a D^{t+1}. \]
The point $\alpha_D$ canonically 
induces a point $\beta_D \in \CX(\bar{k})$, and we have
\[ \deg \beta_D = [k:\Q] \deg \alpha_D, \quad 
   h(\beta_D) = [k:\Q] h(\alpha_D). \]
Finally,
\[ \log |\beta_D,\theta| \leq \log|\alpha_D,\theta_1| \leq - b a D^{t+1}. \] 


\section{Algebraic independence Critieria and Outlook} \label{algind}

As a corollary of Theorem \ref{main}, we get a new proof or the
the Philippon criterion for algebraic independence.

\satz{Theorem} 
For $n \in \N$ let $D_n \in \N,T_n,V_n \in \R^{>0}$ be such that
\[ \lim_{n \to \infty} \frac{V_n}{D_n^k T_n} = \infty, \]
and assume that for each $n$ there is a set global sections 
$\CF_n \subset \Gamma(\Pe^m,O(D_n))$ with
\[ \log |f| \leq T_n, \quad \forall f \in \CF_n, \]
and
\[ \log |\la f|\theta\ra| \leq - V_n, \quad \forall f \in \CF_n, \]
for some $\theta \in \Pe^m(\C)$, and 
\[ \log \left| \bigcap_{f \in \CF_n} V(f),\theta\right| \geq -V_{n-1}, \]
where $V(f´)$ denotes the vanishing set of $f$.
Then the transcendence degree of the
field generated by the coordinates of $\theta$ is at least $k+1$.
\end{Satz}

\vspace{2mm}

{\sc Sketch of Proof:}
Let $\CX$ be the algebraic closure of $\{ \theta \}$, denote by $t$ its
dimension.
For appropriate $K_n$ in relation to $D_n, T_n,V_n$, one uses Theorem
\ref{main} to find a point $\alpha_n$ with 
height and degree bounded in terms of
$K_n$ with logarithmic distance to $\theta$ smaller
than $-V_{n-1}$. Thus, there is an $f \in \CF_n$ such that
$|\alpha_n,\theta|<|\alpha, \di f|$, implying $\alpha_n \notin V(f)$
and Theorem \ref{bezexthe}.2 gives an estimate
\[ 0 = D(\alpha_n . \di f, \theta) \]
against the degree and height of $\alpha_n$, and $f$, and $V_n$ which 
by appropriate choice of $K_n$ can be proved to be less than $0$ if
$t \leq k$.

Actually this proof morover shows that for  
$\theta_1, \cdots, \theta_{k+1}$ any subset of the coordinates of $\theta$
of cardinality $k+1$ with $trdeg_\Q \Q(\theta_1, \ldots, \theta_{k+1})= k+1$
the point $(\theta_1, \ldots, \theta_{k+1})$ forms an $S$-point in
the Zariski closure $\overline{(\theta_1, \ldots, \theta_{k+1})}$ in 
the sense of Mahler classification. For this reason, the details of the
proof will given in \cite{Mahler}.

\vspace{2mm}

There is a stronger version of the main Theorem \ref{main}: Under 
the assumptions of Theorem \ref{main}, assume that $\theta$ is not a
$U$-point in the sense of Mahler classification. Let $a$ be a sufficiently big 
real number and choose a natural number $L$. Then, there is an infinite
subset $M \subset \N$ depending on $a$ and $L$
such that for every $D \in M$ there is a zero dimensional
$\CO_k$-irreducible variety $\CZ_D$ with $d := \deg Z_D \geq L$ such that with
$\alpha_1, \ldots, \alpha_d$ the $\bar{k}$-valued points of $Z_D$, ordered
in such a way that that $|\alpha_1,\theta| \leq \cdots \leq |\alpha_d,\theta|$,
the inqualities
\[ \deg Z_D \leq D^t, \quad h(\CZ_D) \leq a D^t, \]
and
\[ \log |\alpha_1,\theta| \leq \cdots \leq |\alpha_L,\theta| \leq
   - \frac bL \frac{a}{\left( \frac{h(\CX)}a + 2 \deg X \right)^{\frac{t+1}t}}
   \; D^{t+1}, \]
hold, with $b$ a positive real number, again only depending on $t$.
To prove this, one has to consider derivatives of the algebraic distance
which are related by a derivative metric B\'ezout Theorem to be proved in
the next part of this series of papers.

\vspace{2mm}

Theorem \ref{main} allows to prove new algebraic independence criteria, some
of which also depend on the derivative metric B\'ezout Theorem. They also
involve global sections with small evaluations at some point $\theta$ on
a variety, and possibly also small higher derivatives at this point, but
unlike the Philippon criterion above will not require lower bounds on
the distances of the divisors of these global sections to $\theta$.

\begin{appendix}

\section{Proof of Proposition \ref{alg}}

The Proposition will be proved by complete induction on $t$. For 
$t= 1$, and $\CX= \Pe^1$, by \ref{bezout}.1, and the arithmetic B\'ezout
Theorem,
\begin{eqnarray*}
\log |\alpha,\beta| &\geq& D(\alpha,\beta) - c  \deg \alpha \geq
   h(\alpha. \beta) - \deg \alpha h(\beta) - \deg \beta h(\alpha) -
   c' \deg \alpha \deg \beta \\ 
&\geq & - (c' \deg \alpha
   + h(\alpha)) \deg \beta - \deg \alpha h(\beta)- c \deg \alpha, 
\end{eqnarray*}
and the Proposition holds with 
$c_1 = c \deg \alpha +c' \deg \alpha+h(\alpha)$, and $c_2 = \deg \alpha$.

Assume now the Proposition is proved for $t-1$, and consider $\CX = \Pe^t$.
By \cite{App2}, Proposition 4.4, with $D = \sqrt[t]{2 t! \deg \alpha}$
there is a global section $f \in \Gamma(\Pe^t,O(D))$ with
$\log |f| \leq \frac{\sqrt[t]{2 t!}h(\alpha)}{(\deg \alpha)^{\frac{t-a}t}}$.
such that $\di f$ contains $\alpha$. If $\CY:=\di f$ also contains $\beta$,
by induction hypothesis, there are $c_1,c_2>0$ only depending on $t$,
the degree and height of $\alpha$, and the degree and height of $\di f$,
hence by the choice of $f$ only on $t$, and the degree and height of $\alpha$,
such that
\[ \log |\alpha,\beta| \geq -c_1 \deg \beta - c_2 h(\beta). \]
If $\CY$ does not contian $\beta$, then by \ref{bezout}, and the
arithmetic B\'ezout Theorem,
\begin{eqnarray*} 
\log |\alpha,\beta| \geq \log |Y,\beta| &\geq& D(Y,\beta) - c \deg Y \\
&\geq & h(Y.\beta) - \deg Y h(\beta) - h(\CY) \deg \beta - c' \deg Y \deg 
   \beta \\ 
&\geq& - (c' \deg Y + h(\CY)) \deg \beta - \deg Y h(\beta) - c \; deg Y, 
\end{eqnarray*}
and the Proposition follows with $c_1 = c' \deg Y + c \deg Y+ h(\CY)$, 
and $c_2 = \deg Y$, both depending only on $t$, and the degree and
height of $\alpha$.

Let now $\CX$ be arbitrary of relative dimension $t$. 
Embedding $\CX$ into some projective space, and repeating the
argument from the proof of the general case of Theorem \ref{main}, one
obtains the $c_1,c_2$ depending only on $t$, the degree and height
of $\alpha$, and the degree and height of $\CX$ such that
\[ \log |\alpha,\beta| \geq - c_1 \deg \beta - c_2 h(\beta). \]

\end{appendix}

\newpage

\section{Diophantine Approximation on varieties: Leitfaden}

\special{em:linewidth 0.4pt}
\linethickness{0.4pt} 
\vspace{4mm}
\begin{picture}(400,400.00)
\put(70.00,380.00){\makebox{I}}
\put(20,365){\makebox{Algebraic distance and}}
\put(20,350){\makebox{metric B\'ezout Theorem}}
\put(90,335){\line(3,-1){150}}
\put(80,335){\line(1,-6){19}}
\put(325,380){\makebox{II}}
\put(270,365){\makebox{Explicit estimates for ar-}}
\put(270,350){\makebox{ithmetic Hilbert functions}}
\put(320,335){\line(-1,-1){50}}
\put(255,270){\makebox{III}}
\put(190,255){\makebox{Approximation of non-algebraic}}
\put(190,240){\makebox{points by algebraic points}}
\put(260,225){\line(-1,-6){17}}
\put(92,200){\makebox{IV}}
\put(20,185){\makebox{Derivated algebraic distance and}}
\put(20,170){\makebox{derivative metric B\'ezout Theorem}}
\put(110,155){\line(3,-1){110}}
\put(235,100){\makebox{V}}
\put(160,85){\makebox{Algebraic independence criteria}}

\end{picture}

\end{document}